\documentclass{article}
\usepackage{amssymb,amsmath,mathrsfs,latexsym,amscd}
\usepackage{exscale,latexsym,amsthm,graphics}

\usepackage{makeidx}
\makeindex
\usepackage{epsfig}
\usepackage{txfonts}
\newtheorem{thm}{Theorem}[section]
\newtheorem{cor}[thm]{Corollary}

\newtheorem{lemma}[thm]{Lemma}
\newtheorem{prop}[thm]{Proposition}
\newtheorem{proposition}[thm]{Proposition}

\newtheorem{definition}[thm]{Definition}

\newtheorem{remark}[thm]{Remark}

\def\bM {{\mathbb M}}  
 \def\Q {{\mathbb Q}}

\def\N {{\mathbb N}}

\def\pf{ {\medskip \noindent {\bf Proof.}}  \ }
\def\qed{{\hfill $\Box$ \bigskip}}
\def\N {{\mathbb N}}
\def\R {{\mathbb R}}

\def\EE{{\mathbb E}}
\def\P{{\mathbb P}}
\newcommand{\F}{\mathcal{F}}

\def\E{{\mathcal E}}
\def\P{{\mathbb P}}

\def\grad{{\nabla}}

\begin{document}
\noindent
{\Large\bf On the stochastic regularity of distorted Brownian motions}

\bigskip
\noindent
{\bf Jiyong Shin and Gerald Trutnau}
\\

\noindent
{\small{\bf Abstract.} We systematically develop general tools to apply Fukushima's absolute continuity condition. 
These tools comprise methods to obtain a Hunt process on a locally compact separable metric state space  whose transition function has a density w.r.t. the reference measure and methods to estimate drift potentials comfortably. 
We then apply our results to distorted Brownian 
motions and construct weak solutions to singular stochastic differential equations, i.e. equations with possibly unbounded and 
discontinuous drift and reflection terms which may be the sum of countably many local times. The solutions can start from any point of the explicitly specified state space.  
We consider different kind of weights, like Muckenhoupt $A_2$ weights and weights with moderate growth at singularities as well as different kind of (multiple) boundary 
conditions. Our approach leads in particular to the construction and explicit identification of countably skew reflected and normally reflected Brownian motions with singular drift in bounded and unbounded multi-dimensional domains.}\\ 



\section{Introduction}\label{1}
Let $E\subset\R^d$ and $\psi:E\to \R$ be a measurable function such that $\psi>0$ $dx$-a.e. on $E$. We consider a regular Dirichlet form $(\E, D(\E))$ on $L^2(E, \psi dx)$ that can be written as
\begin{equation}\label{i1}
\E (f,g)  = \frac{1}{2}\int_{E}\nabla f\cdot\nabla g \,\psi \,   dx, \quad f,g \in D(\E).
\end{equation} 
The regularity of $(\E, D(\E))$ provides the existence of a Hunt process $\bM = ((X_t)_{t\geq0} , (\P_x)_{x \in E_{\Delta} })$ 
with lifetime $\zeta$ that is associated with $(\E, D(\E))$ and whose generator is informally given as 
$$
Lf=\frac12\Delta f+\frac{\nabla\psi}{2\psi}\cdot \nabla f.
$$ 
$\bM$ is called distorted Brownian motion (cf. \cite{ahkrstr}, \cite{fuku81}, \cite{fuku85}) and forms as in (\ref{1}) with infinitesimal generator $L$ can be generalized to all kind of different 
state spaces $E$ by finding an appropriate interpretation of the gradient $\nabla$ and Laplacian $\Delta$. Due to the good structural properties, like e.g. the self-adjointness 
of the corresponding generators, there is a huge literature about distorted Brownian motion in finite, as well as 
in infinite dimensions (see e.g. \cite{BGG}, \cite{fuku99}, \cite{pos}, \cite{Pan}, \cite{AKR99}, \cite{AnRe} and references therein). 
We shall be concerned with a locally compact separable metric space $E$ for our general results and with $E\subset \R^d$ like above 
in our concrete applications. 
The distorted Brownian motion has then typically an unbounded and discontinuous drift and of special interest is therefore the 
identification of the stochastic differential equation (hereafter SDE) that is fulfilled by it. It is well known how to identify the distorted Brownian motion for quasi-every 
starting point by using Fukushima's decomposition of additive functionals (see \cite{fuku81}, \cite{fuku85}, \cite[Theorem 5.5.1]{FOT}, and \cite{ar91}, \cite[Theorem 2.5]{MR} 
for infinite dimensional state space). 
This approach is in some sense abstract since the set of starting points that is excluded is not explicitly known and rather only given as a set of zero capacity. It can 
nonetheless be made explicit  
by looking at probability distributions $\P_{\nu}(\cdot):=\int_E \P_x(\cdot)\nu(dx)$ where $\nu$ is an explicitly given probability measure that does not charge sets of zero capacity.
Another approach is to solve a corresponding martingale problem for as much as possible explicitly specified starting points (see \cite{AKR}, \cite{BG}, \cite{BGS}, \cite{FaGr}). 
This may be a reasonable intermediate approach, especially if the functions for which the martingale problem is considered are dense in $D(\E)$, but it does not 
lead directly to the identification of the SDE.  
Our strategy for the identification of the distorted Brownian motion for as much as possible explicitly specified starting points is based on Fukushima's absolute continuity 
condition and is known as the strict Fukushima decomposition (cf. \cite[(4.2.9) and Theorem 5.5.5]{FOT},\cite{fuku93}, \cite{fuku94}). To our knowledge it is the first time it is applied systematically 
for weights $\psi\not\equiv const$. For some examples with $\psi\equiv const$, we refer to \cite{BH2}, \cite{fukutomi} and \cite{fuku93}, see also \cite[Examples 5.2.2 and 5.5.3]{FOT}.
The strategy consists of two parts. The first one is to construct a Hunt process whose transition function has a density $p_t(x,y)$ w.r.t. the reference measure $m:=\psi dx$ and is an $m$-version of the $L^2(E,m)$-semigroup $(T_t)_{t>0}$ associated with $(\E, D(\E))$, 
i.e. we need to construct a Hunt process $\bM = (\Omega , \mathcal{F}, (\mathcal{F}_t)_{t\geq0}, (X_t)_{t\geq0} , (\P_x)_{x \in E_{\Delta} })$ with life time $\zeta$ such that
\begin{equation}\label{i2}
P_t f(x) := \mathbb{E}_x[f(X_t)]= \int_{E} p_t(x,y)\, f(y) \, m(dy) 
\end{equation} 
for any $t>0, x \in E$, $f \in \mathcal{B}_b(E)$ and such that $P_t f$ is an $m$-version of $T_t f$ for any $f\in L^2(E,m)\cap \mathcal{B}_b(E)$ and $t>0$.  
Note that even if $(T_t)_{t>0}$ is strong Feller, i.e. $T_t f$ has a continuous $m$-version for any $f \in \mathcal{B}_b(E)$ and $t>0$, so that $T_t f$ has a density as in (\ref{i2}), 
the process constructed via regularity by Dirichlet form methods does not necessarily satisfy 
this condition. In fact since such a Hunt process is only unique for quasi-every starting point (see \cite[Theorem 4.2.8]{FOT}), the absolute continuity condition may be violated 
for some points $x\in E$ in a capacity zero set. 
For the construction of a Hunt process $\bM$ 
on a general locally compact separable metric space $E$ that satisfies the absolute continuity condition, we use two methods. The first one is the well known Feller semigroup 
method that we summarize in Section \ref{2.1.1} and that we apply in the form 
of Lemma \ref{t;Feller}. We then use heat kernel estimates to verify the conditions of Lemma \ref{t;Feller} for concrete Muckenhoupt 
weights (cf. Remark \ref{Fellerapplication}). 
The second method which is developed in Section \ref{2.1.2} is what we call the Dirichlet form method and it is a refinement of the method introduced in 
\cite[Section 4]{AKR}. Our contribution here is to exploit the structure of 
a carr\'e du champ (see Lemma \ref{l;genand} and Remark \ref{normalpro}(ii)) and to find an adequate condition to determine convergence (see ${(\textbf{H2})^{\prime}}$(i) below and proof of Lemma \ref{akrlemma}). For other work, where the method of \cite{AKR} is adopted, we 
refer to \cite{BGS,BG}. As in the case of Feller semigroups, we apply these general results in 
Section \ref{3} to concrete Muckenhoupt $A_2$ weights 
(see Lemma \ref{l;Feller}(i) and Propositions \ref{t;eprocess}, \ref{t;eprocess2}). We remark that it remains open whether the absolute continuity condition holds for general 
Muckenhoupt $A_2$ weights or not. 
According to Proposition \ref{p;strongfl}(i) and (iii), when using the Feller method it remains to show Lemma \ref{t;Feller}(i), and according to 
Proposition \ref{p;strongfl}(i) and (ii), when using the Dirichlet form method it remains to show  ${(\textbf{H2})^{\prime}}$(i) and (ii). In Section \ref{4}, 
we obtain the absolute continuity condition from results of \cite{AKR} using the appropriate part Dirichlet form 
(see Lemma \ref{l;part2}). 
In Section \ref{5}, we assume the absolute continuity condition to be verified, but refer to \cite{BG} to which it accordingly holds under certain conditions 
(see Remark \ref{r;FG}). The results of Section \ref{5}  are also achieved by specifying the appropriate part Dirichlet form (see Lemma \ref{l;part3}). 
The necessary tools for part Dirichlet forms and general auxiliary results are presented in Section \ref{2.2}.
\\ 
The second part of the strategy consists in finding good estimates for the drift potentials 
$$
R_1\mu(x)=\int_E r_1(x,y)\, \mu(dy)
$$
corresponding to the logarithmic derivative 
$\mu:= \frac{\nabla\psi}{2\psi}$ in the sense of distributions and to measures $\mu$ on $\partial E$ that occur through integration 
by parts as boundary terms in case of existing boundary $\partial E$. 
Here $r_1(x,y)=\int_0^{\infty} e^{-t} p_t(x,y)\,dt$. Concretely, in Section \ref{3}, we consider Muckenhoupt $A_2$ weights $\psi=\rho\phi$, where $\rho$ 
is a weakly differentiable function and $\phi$ is a function that is piecewise 
constant and  has discontinuities along boundaries of Euclidean balls (see (\ref{eqphi})), along the boundary of a Lipschitz domain (see (\ref{eqphi2})) 
and along hyperplanes (see (\ref{eqphi3})). 
In this case using informally the Leibniz rule for $\nabla(\rho\phi)$, we see that $\mu(dy)$ is given as the sum of the absolutely continuous part $\frac{\nabla\rho}{2\rho}(y)\,m(dy)$ 
and the corresponding boundary measures. 
In Section \ref{4}, we consider the case where $\phi\equiv 1$ and $E$ has no boundary so that $\mu(dy)=\frac{\nabla\rho}{2\rho}(y)\,m(dy)$ and in Section \ref{5}, 
we consider the case where $\phi\equiv 1$ and 
existing boundary, so that $\mu$ is given as the sum of $\frac{\nabla\rho}{2\rho}(y)\,m(dy)$ and a weighted surface measure (see Lemma \ref{ibp3}). 
Our key for estimating potentials is Proposition \ref{p;smooth} that we found very useful and apply throughout the article. Especially, if no continuity properties 
of a potential are known, we use resolvent kernel estimates to find continuous Riesz potentials (see (\ref{rieszpot}) and Lemma \ref{l;miz}) as upper bound $r_1^G$ as in  
Proposition \ref{p;smooth} for the potential, i.e. we use Proposition \ref{p;smooth} in combination with resolvent kernel estimates and Lemma \ref{l;miz}.
We use this procedure for instance globally in Lemma \ref{l;Feller}(iii)-(v) where for the global resolvent kernel estimates, we use known global heat kernel estimates 
for Muckenhoupt weights from \cite{St3} (see (\ref{heatesti})). We use it locally in Lemma \ref{t;sm3} using
local heat kernel estimates that we derive using Nash type inequalities and the Davies method of \cite{CKS} similarly to what is done in \cite[Theorems 2.3, 3.1]{BH2}
(see Lemma \ref{T;NI1}, Proposition \ref{t;tde1} and Corollary \ref{C;RDE1}). 
Of special interest could be the corresponding localization procedure via part processes that we apply on a nice exhaustive sequence of sets for the state space 
(see conditions ($\iota$), ($\kappa$) in Section \ref{5}, Lemma \ref{t;sm3}, Proposition \ref{t;lsfd3}, Lemmas \ref{l;limit}, \ref{T;LOCAL} and proof of Theorem \ref{t;lsfd4}). 
We use it when global resolvent kernel estimates do not provide enough regularity or are not at hand. For other places in this article where we use this localization procedure 
see Proposition \ref{p;3.8}(ii), Theorem \ref{t;3.9}(ii) and Remark \ref{r;3.15}.\\ 
The Muckenhoupt $A_2$ weights $\psi=\rho\phi$ that we investigate in Section \ref{3}, lead to solutions of SDEs of the following type
\begin{equation}\label{dbm}
X_{t} = x + W_{t} + \int_{0}^{t}\frac{\nabla \rho}{2 \,\rho}(X_s) \, ds + L_t^{\phi},\ \  t  \ge 0, \ x\in E\subset \R^d,
\end{equation}
where $L^{\phi}$ may be a series of local times (see (\ref{LSFD1}) of Theorem \ref{t;Ldec}). Theorem \ref{t;Ldec} is formulated under general conditions on $\rho$ and $\phi$. 
We then extensively study  the typical case of an $A_2$ weight where $\rho(x)=\|x\|^{\alpha}$, $\alpha\in (-d+1,d)$ and $\phi$ is an explicitly given piecewise constant function that is globally 
bounded above and below by strictly positive constants. 
In this case it is known that the capacity of $\{0\}$ is zero, iff $\alpha\in [-d+2,d)$. We obtain that one can choose 
$E=\R^d$, if $\alpha\in (-d+1,2)$ and $L^{\phi}\equiv 0$,  or if $\alpha\in (-d+1,1)$ and $L^{\phi}\not\equiv 0$ (see Proposition \ref{p;3.8}(i), Theorem \ref{t;3.9}(i) and Theorem \ref{t;Sdec}) and that 
one can choose $E=\R^d\setminus\{0\}$ in the remaining cases (see Proposition \ref{p;3.8}(ii), Theorem \ref{t;3.9}(ii) and Remark \ref{r;3.15}). 
Two observations are here worth to be noted. The first is that we are able to start 
in $0$ although $\{0\}$ might be a capacity zero set and the second is that we lose one dimension in $\alpha$ in case there are boundary terms. The reason for the last is that we use continuous 
Riesz potentials of the form (\ref{rieszpot}) as upper bounds for our drift potentials and that drifts which are given as surface measures on a nice boundary are 
equivalent to the Lebesgue measure of one dimension less (cf. Lemma \ref{l;Feller}(v)). The concrete examples of drifts $L^{\phi}$ that we obtain in (\ref{dbm}) 
can be summarized as follows. If $\phi$ is as in (\ref{eqphi}) piecewise constant on countably many annuli with jumps along their boundaries, $L^{\phi}$ is given as the last term 
in (\ref{equast}) which corresponds to a distorted Brownian motion with skew reflection on the boundary of Euclidean balls that may accumulate. (\ref{equast}) seems new to us. We could not find any similar equation in multi-dimensions in the literature. Its one-dimensional counterpart is studied extensively in \cite{OuRuTr13}.
If $\phi$ is as in (\ref{eqphi2}) piecewise constant on a bounded Lipschitz domain and on its complement, then $L^{\phi}$ is 
given as a scalar multiple of the boundary local time on the boundary of a Lipschitz domain $G$ as in (\ref{SFD1}). The corresponding process could be called a $\beta$-skew distorted Brownian motion w.r.t. $G$.
In case of skew reflection at the boundary of a $C^{1,\lambda}$-domain, $\lambda\in(0,1]$ and smooth diffusion coefficient, a weak solution is constructed in 
\cite[III. \S 3 and \S 4]{port}, see also references therein. The reflection term in \cite{port} is defined as generalized drift. 
If $\phi$ is as in (\ref{eqphi3}) piecewise constant on countably many infinite strips with jumps along countably many hyperplanes, then $L^{\phi}$ is given as the last term in (\ref{LSFD2}). 
Variants of (\ref{LSFD2}), but without accumulation points and Lipschitz drift appear in \cite{Tom, O82, taka86a}. For recent related work, we refer to \cite{atbu}.\\
In Section \ref{4}, we complete results of \cite{AKR}. There the distorted Brownian motion is constructed on $\R^d\setminus \{\psi=0\}$ for certain weights $\psi$, 
but the corresponding SDE is not explicitly identified. It was noted in \cite[Remark 5.6]{AKR} that besides using direct stochastic calculus one could possibly also achieve this identification by refining arguments from \cite{FOT}. 
As already mentioned, we work out the latter using the part Dirichlet form on $\R^d\setminus \{\psi=0\}$. For details we refer to Section \ref{4}.\\
In Section \ref{5}, we complete results from \cite{Tr1}. Precisely, under the assumptions $(\eta)-(\kappa)$ of Section \ref{5}, we show in Theorem \ref{t;lsfd4} that the Skorokhod decomposition 
that was obtained in \cite{Tr1} for quasi-every starting point can be achieved in concrete examples for every starting point outside an explicitly specified capacity zero set in the symmetric case. We note that the absolute continuity condition is assumed to hold  in $(\theta)$.
For additional conditions according to which the absolute continuity condition is satisfied, we refer to \cite{BG} (see Remark \ref{r;FG}). For work that is strongly related with 
Theorem \ref{t;lsfd4}, we refer to \cite{BH2, ch, fukutomi, PaWi}.\\
Finally, let us remark that we only treat the semimartingale case, but that the strict Fukushima decomposition has also been formulated in the non-semimartingale case (see \cite{fuku93}). 
It could be interesting to see which phenomena occur in this case. Moreover, because we did not 
want to overload this presentation, we also did not consider the $(a_{ij})$-case in our concrete examples. 
But drift potentials that occur in the $(a_{ij})$-case can be handled by exactly the same methods 
that are presented here once the absolute continuity condition is established. For this, we refer to forthcoming work.
\section{Preliminaries and the absolute continuity condition}\label{2}
Let $d\ge1$. $C_0^{\infty}(\R^d)$ denotes the set of all infinitely differentiable functions with compact support in $\R^d$. Let $\nabla f : = ( \partial_{1} f, \dots , \partial_{d} f )$  and  $\Delta f : = \sum_{j=1}^{d} \partial_{jj} f$ where $\partial_j f$ is the $j$-th weak partial derivative of $f$ and $\partial_{jj} f := \partial_{j}(\partial_{j} f) $, $j=1, \dots, d$.
As usual $dx$ is the Lebesgue measure on $\R^d$ and $\delta_x$ is the Dirac measure at $x$.  For any open set $G \subset \R^d$ the Sobolev space $H^{1,q}(G, dx)$, $q \ge 1$ is defined to be the set of all functions $f \in L^{q}(G, dx)$ such that $\partial_{j} f \in L^{q}(G, dx)$, $j=1, \dots, d$, and $H^{1,q}_{loc}(\R^d, dx) : =  \{ f  \,|\;  f \cdot 1_{U} \in H^{1,q}(U, dx),\,\forall U \subset \R^d,\, U \; \text{relatively}$ $\; \text{compact} \; \text{open}  \}$. We always equip $\R^d$ with the Euclidean norm $\| \cdot \|$ and write $B_{r}: = \{ x \in \R^d \ | \ \|x\| < r  \}$.   \\
 
For a locally compact separable metric space $(E,d)$ with Borel  $\sigma$-algebra $\mathcal{B}(E)$ we denote the set of all $\mathcal{B}(E)$-measurable $f : E \rightarrow \R$ which are bounded, or nonnegative by $\mathcal{B}_b(E)$, $\mathcal{B}^{+}(E)$ respectively. $B_{r}(y): = \{ x \in E \ | \ d(x,y) < r  \}$, $r>0$, $y \in E$. $L^q(E, \mu)$, $q \in[1,\infty]$ are the usual $L^q$-spaces equipped with $L^{q}$-norm $\| \cdot \|_{q}$ with respect to the  measure $\mu$ on $E$, $\mathcal{A}_b$ : = $\mathcal{A} \cap \mathcal{B}_b(E)$ for $\mathcal{A} \subset L^q(E,\mu)$,  and $L^{q}_{loc}(E,\mu) := \{ f \,|\; f \cdot 1_{U} \in L^q(E, \mu),\,\forall U \subset E, U \text{ relatively compact open} \}$, where $1_A$ denotes the indicator function of a set $A$. As usual, we also denote the set of continuous functions on $E$, the set of continuous bounded functions on $E$, the set of compactly supported continuous functions in $E$ by $C(E)$, $C_b(E)$, $C_0(E)$, respectively. $C_{\infty}(E)$ denotes the space of continuous functions on $E$ which vanish at infinity. For $A \subset E$ let $\overline{A}$ denote the closure of $A$ in $E$, $A^c:=E\setminus A$. We will refer to \cite{FOT} till the end, hence some of its standard notations may be adopted below without definition.\\

 In order to simplify notation while handling inequalities or estimates we make the convention that unless otherwise specified $c>0$ stands for an arbitrary constant whose value may vary from inequality to inequality.
\subsection{Global setting}
Throughout, we let $(\E,D(\E))$ be a symmetric, strongly local, regular Dirichlet form on $L^{2}(E, m)$ where $m$ is a positive Radon measure on $(E, \mathcal{B}(E))$ with full support on $E$. We further assume throughout that $\E$ admits a carr\'e du champ
\[
\Gamma: D(\E) \times D(\E) \to L^1(E,m)
\]
as in \cite[Definition 4.1.2]{BH}. As usual we define $\E_1(f,g) := \E(f,g) + (f,g)_{L^{2}(E , \, m)}$ for $f,g \in D(\E)$ and  $\| \,f\, \|_{D(\mathcal{E})} : = \E_1(f,f)^{1/2},  \;   f \in D(\E)$.
Let  $(T_t)_{t > 0}$ and $(G_{\alpha})_{\alpha > 0}$ be the $L^2(E, m)$-semigroup and resolvent associated to $(\E,D(\E))$ and $(L,D(L))$ be the corresponding generator (see \cite[Diagram 3, p. 39]{MR}). Let Cap be the capacity related to the regular symmetric Dirichlet form $(\E,D(\E))$ as defined in \cite[2.1]{FOT}. We say that a function $f$ is locally in $D(\E)_b$ ($f \in D(\E)_{b,loc}$ in notation) if for any relatively compact open set $G \subset E$, there exists a function $g \in D(\E)_b$
such that $f=g \;\; m$-a.e. on $G$. We consider the condition \\

$(\bf{H1})$ There exists a $\mathcal{B}(E) \times \mathcal{B}(E)$ measurable non-negative map $p_{t}(x,y)$ such that
\[
P_t f(x) := \int_{E} p_t(x,y)\, f(y) \, m(dy) \,, \; t>0, \ \ x \in E,  \ \ f \in \mathcal{B}_b(E),
\]
is a (temporally homogeneous) sub-Markovian transition function (see \cite[Section 1.2]{CW}) and an $m$-version of $T_t f$ if $f  \in  L^2(E , m)_b$.\\
\newline
$p_{t}(x,y)$ is called the transition kernel density or heat kernel. Taking the Laplace transform of $p_{\cdot}(x, y)$, we see that $(\bf{H1})$ implies that there exists a $\mathcal{B}(E) \times \mathcal{B}(E)$ measurable non-negative map $r_{\alpha}(x,y)$ such that
\begin{equation}\label{Rdensity}
R_{\alpha} f(x) := \int_{E} r_{\alpha}(x,y)\, f(y) \, m(dy) \,, \; \alpha>0, \; x \in E, f \in \mathcal{B}_b(E),
\end{equation}
is an $m$-version of $G_{\alpha} f$ if $f  \in  L^2(E, m)_b$. $r_{\alpha}(x,y)$ is called the resolvent kernel density.
For a signed Radon measure $\mu$ on $E$, let us define
\begin{equation}\label{Rpoten}
R_{\alpha} \mu (x) = \int_{E} r_{\alpha}(x,y) \, \mu(dy) \, , \;\; \alpha>0, \;\;x \in E,
\end{equation}
whenever this makes sense. Throughout, we set $P_0 : = id$.
Furthermore, assuming that $(\bf{H1})$ holds, we can consider the condition\\ 

$(\bf{H2})$ There exists a Hunt process with transition function $(P_t)_{t \ge 0}$.\\
\newline
We recall that $(\bf{H2})$ means that there exists a Hunt process
\begin{equation}
\bM = (\Omega , \mathcal{F}, (\mathcal{F}_t)_{t\geq0}, \zeta ,(X_t)_{t\geq0} , (\P_x)_{x \in E_{\Delta} }),
\end{equation}
with state space $E$ and life time $\zeta$ such that $P_t(x,B) : = P_t 1_B (x) = \P_x(X_t \in B)$ for any $x \in E$, $B \in \mathcal{B}(E)$, $t \ge 0$. Here, $\Delta$ is the cemetery point and as usual any function $f : E \rightarrow \R$ is extended to $\{\Delta\}$ by setting $f(\Delta):=0$. $E_{\Delta}: = E \cup \{\Delta\}$ is the one-point compactification if $E$ is not already compact, if $E$ is compact then $\Delta$ is added to $E$ as an isolated point. 

\begin{remark}\label{r;absolute}
Note that if $(\bf{H1})$ and $(\bf{H2})$ hold, then $\bM$ is associated with $(\E, D(\E))$ and satisfies the {\it absolute continuity condition} as stated in \cite[p. 165]{FOT}.
\end{remark}

Below, we present two methods to obtain $\bM$ as in Remark \ref{r;absolute}.

\subsubsection{The Feller semigroup method}\label{2.1.1}

Assuming (\textbf{H1}), a Hunt process as in ($\bf{H2}$) can be constructed by means of a Feller semigroup (cf. \cite[(9.4) Theorem]{BlGe}). For the definition of Feller semigroup, we refer to \cite[Section 2.2]{CW}.
\begin{remark}\label{r;feller}
Under $(\bf{H1})$, $(P_t)_{t \ge 0}$ is a Feller semigroup, if 
\begin{itemize}
\item[(i)] $\forall f \in C_{\infty}(E)$, $\lim_{t \to 0} P_t f= f$ uniformly on $E$,
\item[(ii)] $P_t C_{\infty}(E) \subset C_{\infty}(E)$ for each $t > 0$.
\end{itemize}
\end{remark}
It is well known that the condition of uniform convergence in Remark \ref{r;feller} (i) can be relaxed to pointwise convergence (see for instance \cite[Section 2.2 Exercise 4.]{CW}). 
The conditions of Remark \ref{r;feller} can be further relaxed to the conditions of the following lemma which are suitable for us. 
\begin{lemma}\label{t;Feller}
Suppose $(\bf{H1})$ and that
\begin{itemize}
\item[(i)] $\lim_{t \to 0} P_t f(x) = f(x)$ for each $x \in E$ and $f \in C_0(E)$,
\item[(ii)] $P_t C_0(E) \subset C_{\infty}(E)$ for each $t>0$.
\end{itemize} 
Then $(P_t)_{t \ge 0}$ is a Feller semigroup. In particular (\textbf{H2}) holds.
\end{lemma}

\begin{remark}\label{Fellerapplication}
One can use heat kernel estimates for $p_t(x,y)$ to check the assumptions of Lemma \ref{t;Feller} (i), (ii) (see Lemma \ref{l;Feller} (i) below).
\end{remark}

\subsubsection{The Dirichlet form method}\label{2.1.2}
The second method to obtain a Hunt process as in Remark \ref{r;absolute}, given a transition function as in ($\textbf{H1}$), is by a method that we shall call the Dirichlet form method.
It is a refinement of the method introduced in \cite[Section 4]{AKR}. We shall put it in a frame that is suitable for our purposes. We assume hence ($\textbf{H1}$) to hold and explain the main steps of the method and of our refinement.
\newline 
Given the transition function $(P_t)_{t \ge 0}$ on $E$, restricted to the positive dyadic rationals $S:= \bigcup_{n \in \N} S_n$, $S_n : = \{k2^{-n} \ | \ k \in \N \cup \{0\}\}$, we construct a Markov process 
\[
\bM^{0} = (\Omega , \mathcal{F}^0, (\mathcal{F}^0_s)_{s \in S}, (X_s^0)_{s \in S} , (\P_x)_{x \in E_{\Delta} })
\]
with transition function on $E_{\Delta}$
\[
P_t^{\Delta}(x, dy) =
\begin{cases}
\big[1- P_t(x,E)\big] \delta_{\Delta} (dy) + P_t(x, dy), \quad \text{if} \ x \in E\\
\delta_{\Delta} (dy), \quad \text{if} \ x = \Delta
\end{cases}
\]
by Kolmogorov's  method (see \cite[Chapter III]{RYor}).
Here $\Omega : = (E_{\triangle})^S$ is equipped with the product $\sigma$-field $\F^0$, $X_s^0 : (E_{\triangle})^S \to E_{\triangle}$ are coordinate maps and $\F_s^0 : = \sigma( X_r^0 \ | \ r \in S, r \le s)$. 
By the theory of Dirichlet forms there exists a Hunt process
\[
\tilde{\bM} = (\tilde{\Omega
}, \tilde{\mathcal{F}}, (\tilde{\mathcal{F}}_t)_{t\ge 0}, \tilde{\zeta}, (\tilde{X}_t)_{t \ge 0} , (\tilde{\P}_x)_{x \in E_{\Delta} })
\]
associated with $(\mathcal{E},D(\mathcal{E}))$, where $\tilde{\Omega} = \{\omega = (\omega(t))_{t \ge 0} \in C([0,\infty),E_{\Delta}) \ | \ \omega(t) = \Delta, \ \ \forall t \ge \tilde{\zeta} \}$ (see \cite[Theorem 4.5.3]{FOT}). Let $\nu : = gdm$, where $g \in L^1(E,m)$, $g>0$ $m$-a.e., $\int_E g \ dm =1$, and set
\[
\tilde{\P}_{\nu}(\cdot) := \int_{E} \tilde{\P}_x (\cdot) \ g(x) \ m(dx).
\]
Consider the one-to-one map $G : \tilde{\Omega} \to \Omega$ defined by
\[
G(\omega) = \omega |_S.
\]
Then $G$ is $\tilde{\F}^0 / \F^0$ measurable and $\tilde{\Omega} \in \tilde{\F}^0$, where $\tilde{\F}^0 : = \sigma(\tilde{X}_s \ | \ s \in S)$ and exactly as in \cite[Lemma 4.2 and 4.3]{AKR} we can show that $\tilde{\P}_{\nu} |_{\tilde{\F}^0} \circ G^{-1} = \P_{\nu}$, $G(\tilde{\Omega}) \in \F^0$ and $\P_{\nu}(G(\tilde{\Omega}))=1$. Then, we show Lemma 4.4 of \cite{AKR} with $A = G(\tilde{\Omega}) \ \ \forall x \in E$, i.e. if
\[
\Omega_1 : = \bigcap_{s > 0, s \in S} \theta_s^{-1} (G(\tilde{\Omega})),
\]
where $\theta_s : \Omega \to \Omega$, $\theta_s(\omega) : = \omega(\cdot +s)$, for $s \in S$, is the usual shift operator, then 
\begin{equation}\label{omega1}
\P_x(\Omega_1) = 1
\end{equation}
for all $x \in E$.\\ 
Before we go on with our refinement of the Dirichlet form method it is convenient to introduce some definitions and lemmas:\\
If $\mathcal{A}$ is a set of functions $f : E \to \R$, we define $\mathcal{A}_0 : = \{f \in \mathcal{A} \ | $ supp($f$) : = supp($|f| dm$) is compact in $E \}$. It is well-known that $T_t, \ t>0$, restricted to $L^1(E,m) \cap L^{\infty}(E,m)$ can be extended to a $C_0$-semigroup of sub-Markovian contractions on $L^r(E,m)$ for any $r \ge 1$. We denote the corresponding generators by $(L_r, D(L_r))$ (for details we refer to \cite[Lemmas 1.11 and 1.12 of Appendix B]{Eb} and references therein).
\begin{lemma}\label{l;genand}
Let $u \in D(L)_0 \cap \mathcal{B}_b(E)$. Then:
\begin{itemize}
\item[(i)] \emph{supp}($Lu$) $\subset$ \emph{supp}($u$).
\item[(ii)] It holds $u,u^2 \in D(L_1)$ and
\[
L_1 u^2 = \Gamma(u,u) + 2 u Lu.
\]
\item[(iii)] If $\Gamma(u,u) \in L^p(E,m)$ for some $p \in [2,\infty]$, then $u^2 \in D(L)_0 \cap \mathcal{B}_b(E)$.
\end{itemize}
\end{lemma}
\pf
(i) The statement follows easily from the local property of $(\mathcal{E},D(\mathcal{E}))$, since
\[
\int Lu \cdot \text{v} \ dm = - \mathcal{E}(u,\text{v}) = 0 \quad \forall \text{v} \in D(\mathcal{E}) \ \text{with supp}(\text{v}) \subset \R^d \setminus \text{supp}(u).
\]
(ii) Since $L^2(E,m)_0 \subset L^1(E,m)_0$, we conclude with the help of (i) that $u,Lu \in L^1(E,m)_0$. Hence $u \in D(L_1) \cap \mathcal{B}_b(E)$ by \cite[Lemmas 1.11, 1.12 of Appendix B]{Eb}. By \cite[I. Theorem 4.2.1]{BH}, it then holds $u^2 \in D(L_1) \cap \mathcal{B}_b(E)$ and
\[
L_1 u^2 = \Gamma(u,u) + 2u Lu.
\]
(iii) By \cite[Lemma 3.8 (iii)]{Tr2} we find supp( $\Gamma(u,u))\subset$ supp($u$) since $1_{\R^d \setminus \text{supp}(u)} \Gamma(u,u) dm = 0$. Therefore $ \Gamma(u,u) \in L^2(E,m)_0$ and so $L_1 u^2 \in L^2(E,m)$ by $(ii)$. Since $u^2 \in L^2(E,m)$ and $u^2 \in D(L_1)$ by (ii) it follows again from \cite[Lemmas 1.11, 1.12 of Appendix B]{Eb} that $u^2 \in D(L)$.
\qed
\begin{lemma}\label{l;lemma0}
Let $f \in \mathcal{B}(E)$ such that $R_1|f|$ is finite on $E$ (for instance if $R_1|f|$ is continuous or if $f \in L^{\infty}(E,m)$). Then for any $t \ge 0$
\[
\lim_{\begin{subarray} \ s \downarrow t \\ s \in S  \end{subarray}  } P_s R_1 f (x) = e^t \int_t^\infty e^{-u} \ P_{u} f(x) \ du = P_t R_1 f(x).
\]
In particular
\[
\lim_{\begin{subarray} \ s \downarrow 0 \\ s \in S  \end{subarray}  } P_s R_1 f (x) = R_1 f(x) \quad \text{for any} \ \ x \in E.
\]
\end{lemma}
\pf
First note that for any function $f \in \mathcal{B}^+(E)$, we have $P_s f(x) = P_s^{\Delta} f(x)$ if $x \in E$.
Using this, for any $f \in \mathcal{B}^+(E)$ and $x \in E$, we then obtain with Fubini
\begin{equation}\label{pr1f}
P_s R_1 f(x) = P_s^{\Delta} R_1 f(x) = \EE_x [R_1 f(X_s^0)] = e^s \int_s^{\infty} e^{-u} \ P_u f(x) \ du, \quad s>0,
 \end{equation}
where $\EE_x$ denotes the expectation w.r.t. $\P_x$. The r.h.s. of \eqref{pr1f} converges in $\R$ to $e^t \int_t^{\infty} e^{-u} \ P_u f(x) \ du$ as $s \downarrow t$, $t \ge 0$ if $R_1 f(x)$ is finite. If $R_1 |f|$ is finite, then $R_1(f^+)$ as well as $R_1 (f^-)$ are finite and so the assertion follows.
\qed

$\Omega_1$ defined in \eqref{omega1} consists of paths in $\Omega$ which have unique continuous extensions to $(0,\infty)$ which still lie in $E_{\Delta}$ and stay in $\Delta$ once they have hit $\Delta$. Following the main idea of \cite{AKR}, we have to handle the limits at $s=0$. This can be done assuming the following condition\\

${(\textbf{H2})^{\prime}}$ We can find $\{ u_n \ | \ n \ge 1 \} \subset D(L) \cap C_0(E)$ satisfying:
\begin{itemize}
\item[(i)] For all $\varepsilon \in \Q \cap (0,1)$ and
$y \in D$, where $D$ is any given countable dense set in $E$, there exists $n \in \N$ such that $u_n (z) \ge 1$, for all $z \in \overline{B}_{\frac{\varepsilon}{4}}(y)$ and $u_n \equiv 0$ on $E \setminus B_{\frac{\varepsilon}{2}}(y)$.
\item[(ii)] $R_1\big( [(1 -L) u_n]^+ \big)$, $R_1\big( [(1 -L) u_n]^- \big)$, $R_1 \big( [(1-L_1)u_n^2]^+ \big)$, $R_1 \big( [(1-L_1)u_n^2]^- \big)$ are continuous on $E$ for all $n \ge 1$.
\item[(iii)] $R_1 C_0(E) \subset C(E)$.
\item[(iv)] For any $f \in C_0(E)$ and $x \in E$, the map $t \mapsto P_t f(x)$ is right-continuous on $(0,\infty)$.
\end{itemize}

\begin{remark}\label{normalpro}
\begin{itemize}
\item[(i)] By Lemma \ref{l;genand} (ii), $u_n^2 \in D(L_1)$ $\forall n \ge 1$. Thus $L_1 u_n^2$ in $(\textbf{H2})^{\prime}$ (ii) is well-defined.
\item[(ii)] In view of Lemma \ref{l;genand}  $(\textbf{H2})^{\prime}$ (ii)-(iii) can be replaced by the following (stronger) condition:\\

$\exists r \in [1,\infty]$ such that $R_1 \big( L^r(E,m)_0 \big) \subset C(E)$ and $Lu_n \in L^r(E,m)$ for any $n \ge 1$ and if $r \neq 1$, then $\Gamma(u_n,u_n)^{1/2}  \in L^{\infty}(E,m), \ \forall n \ge 1$.\\
\end{itemize}
\end{remark}

Define
\[
\Omega_0 : = \{ \omega \in \Omega_1 \ | \ \lim_{s \downarrow 0} X_s^0(\omega) \ \text{exists in } \ E\}.
\]

\begin{lemma}\label{akrlemma}
Under $(\textbf{H1})$ and $\textbf{(H2)}^{\prime}$, we have 
\begin{equation}\label{normal}
\lim_{\begin{subarray}{1} s \downarrow 0 \\ s \in S \end{subarray} } X_s^0 =x \quad \P_x\text{-a.s.} \quad \text{for all } \  x \in E.
\end{equation}  
In particular $\P_x (\Omega_0) = 1$ for any $x \in E$.
\end{lemma}
\pf
Let $x \in E$, $n \ge 1$. Then the processes
\[
\Big( e^{-s} R_1 \big( [(1-L) u_n]^+ \big) (X_s^0), \F_s^0, \P_x  \Big) \quad \text{and} \quad \Big( e^{-s} R_1 \big( [(1-L) u_n]^- \big) (X_s^0), \F_s^0, \P_x  \Big)
\]
are positive supermartingales. Indeed since $R_1 \big( [(1-L) u_n]^{\pm}\big)$ is continuous by $\textbf{(H2)}^{\prime}$ (ii), the processes are adapted and integrable. The supermartingale property follows by standard manipulations using the simple Markov property. Then by \cite[1.4 Theorem 1]{CW} for any $t \ge 0$
\[
\exists \lim_{\begin{subarray} \ s \downarrow t \\ s \in S  \end{subarray}} e^{-s} \  R_1 \big( [(1-L) u_n]^{\pm}\big) (X_s^0) \quad \P_x\text{-a.s.}
\]
thus
\begin{equation}\label{existslim}
\exists \lim_{\begin{subarray} \ s \downarrow 0 \\ s \in S \end{subarray}} u_n(X_s^0) \quad \P_x\text{-a.s.}
\end{equation}
We have $u_n = R_1 \big((1-L)u_n \big)$ and $u_n^2 = R_1 \big((1-L_1)u_n^2 \big)$ $m$-a.e., but since both sides are respectively continuous by $\textbf{(H2)}^{\prime}$ (ii) and $m$ has full support, it follows that the equalities hold pointwise on $E$. Therefore
\[
\EE_x \big[\big(u_n(X_s^0) - u_n(x)   \big)^2 \big] = P_s R_1\big( (1-L_1)u_n^2 \big)(x) - 2u_n(x) P_s R_1\big((1-L)u_n\big)(x) + u_n^2(x)
\] 
and so
\begin{equation}\label{separ}
\lim_{\begin{subarray} \ s \downarrow 0 \\ s \in S  \end{subarray}} \EE_x\big[ \big(u_n(X_s^0) -u_n(x) \big)^2 \big] =0
\end{equation}
by Lemma \ref{l;lemma0}.  \eqref{existslim} and \eqref{separ} now imply that 
\begin{equation}\label{separating}
\lim_{\begin{subarray} \ s \downarrow 0 \\ s \in S  \end{subarray}} u_n(X_s^0 (\omega)) = u_n(x) \quad \text{for all} \ \omega \in \Omega_x^n,
\end{equation}
where $\Omega_x^n \subset \Omega_1$ with $\P_x (\Omega_x^n) = 1$. Let $\omega \in \Omega_x^0 : = \bigcap_{n \ge 1} \Omega_x^n$. Then $\P_x (\Omega_x^0) = 1$. Suppose that $X_s^0(\omega)$ does not converge to $x$ as $s \downarrow 0$, $s \in S$. Then there exists $\varepsilon_0 \in \Q$ and a subsequence $(X_{s_k}^0(\omega))_{k \in \N}$ such that $d(X_{s_k}^0(\omega), x) > \varepsilon_0$ for all $k \in \N$. For $\varepsilon_0 \in \Q$ we can find $y \in D$ and $u_n$ in $\textbf{(H2)}^{\prime}$ (i) such that $d(x,y) \le \frac{\varepsilon_0}{4}$ and $u_n(z) \ge 1 $, $z \in  \overline{B}_{\frac{\varepsilon_0}{4}}(y)$  and $u_n(z) = 0$, $z \in E \setminus B_{\frac{\varepsilon_0}{2}}(y)$. Then $u_n(X_{s_k}^0(\omega))$ can not converge to $u_n(x)$ as $k \to \infty$. This is a contradiction.
\qed\\
Now we define for $t \ge 0$
\[
X_t(\omega) : = 
\begin{cases}
\lim_{\begin{subarray}{1} s \downarrow t \\  s \in S \end{subarray} } X_s^0(\omega) \quad \text{if} \ \omega \in \Omega_0 \\
x_0 \quad \text{if} \ \omega \in \Omega \setminus \Omega_0,
\end{cases}
\]
where $x_0$ is an arbitrary but fixed point in $E$. Then by $\textbf{(H2)}^{\prime}$ (iv) for any $t \ge 0$, $f \in C_0(E)$ and $x \in E$ 
\[
\EE_x[f(X_t)] = P_t f(x).
\]
Since $\sigma(C_0(E)) = \mathcal{B}(E)$, it follows that 
\[
\bM = (\Omega, \mathcal{F}, (\mathcal{F}_t)_{t\ge 0}, (X_t)_{t\geq0} , (\P_x)_{x \in E_{\Delta} }),
\]
where $(\mathcal{F}_t)_{t\ge 0}$ is the natural filtration, is a normal Markov process with transition function $(P_t)_{t \ge 0}$. Moreover, $\bM$ has continuous paths up to infinity on $E_{\Delta}$. The strong Markov property of $\bM$ follows from \cite[Section I. Theorem (8.11)]{BlGe} using $\textbf{(H2)}^{\prime}$ (iii). Hence $\bM$ is a Hunt process, i.e. a strong Markov process with continuous sample paths on $E_{\Delta}$, and has $(P_t)_{t \ge 0}$ as transition function. Therefore \textbf{(H2)} holds. Making a statement out of the last conclusion we put it in the following lemma.
\begin{lemma}\label{l;2.10}
Assume \textbf{(H1)} holds. Then $\textbf{(H2)}^{\prime}$ implies \textbf{(H2)}.
\end{lemma}

\begin{remark}
If $(T_t)_{t \ge 0}$ is strong Feller, i.e. $T_t f$ has a continuous $m$-version for any $f \in \mathcal{B}_b(E)$ and $\textbf{(H2)}^{\prime}$ (i)-(ii) and $\textbf{(H2)}^{\prime}$ (iv)
hold, then (\textbf{H1}) and (\textbf{H2}) hold (cf. Proof of Proposition \ref{p;strongfl} below).
\end{remark}

\subsection{Local setting and general auxiliary results}\label{2.2}
We assume ($\bf{H1}$) and ($\bf{H2}$) throughout the Section \ref{2.2}.
\begin{definition}
Let $B$ be an open set in $E$. For $x \in B, t \ge 0, \alpha>0$ and $p \in [1,\infty)$  let
\begin{itemize}
\item $\sigma_{B^c} := \inf\{t>0 \,|\;X_t \in B^c\}$, $D_{B^c} := \inf\{t\ge0 \,|\;X_t \in B^c\}$,
\item $P^{B}_{t}f(x) : = \EE_x [f(X_t) ; t<\sigma_{B^c} ]\; , \; f \in \mathcal{B}_{b}(B)$,
\item $R^{B}_{\alpha}f(x) : = \EE_x \Big[\int_{0}^{\sigma_{B^c}} e^{-\alpha s} f(X_s) \, ds \Big]  \; , \; f \in \mathcal{B}_{b}(B)$ ,
\item $D(\E^{B}): = \{u \in D(\E) \, | \; u=0 \,\,  \mathcal{E}\text{-q.e} \; on  \; B^c \}$.
\item $\E^{B} : = \mathcal{E} \, |_{D(\E^{B})\times D(\E^{B})}$.
\item $L^2(B \, ,m): = \{u \in L^2(\R^d , m) \,|\; u=0, \; m \text{-a.e. on} \; B^c\}$.
\item $||f||^p_{p,B}: = \int_{B} |f|^p \; dm$. 
\item $||f||_{\infty,B}:= \inf \Big\{c>0 \,|\; \int_{B}  1_{\{ \,|f|>c \, \} } \ dm = 0 \Big\}$.
\item $\E^{B}_1(f,g) : = \E^{B}(f,g) + \int_{B} f g \; dm,  \;\;     f,g \in D(\E^{B})$.
\item $\| \,f\, \|_{D(\mathcal{E}^{B})} : = \E^{B}_1(f,f)^{1/2},  \;\;     f \in D(\E^{B})$.
\end{itemize}
\end{definition}
$(\E^{B},D(\E^{B}))$ is called the part Dirichlet form of $(\E, D(\E))$ on $B$. It is a regular Dirichlet form on $L^2(B, m)$ (cf. \cite[Section 4.4]{FOT}).
Let  $(T^{B}_t)_{t > 0}$ and $(G^{B}_{\alpha})_{\alpha > 0}$ be the $L^2(B, m)$-semigroup and resolvent associated to $(\E^{B},D(\E^{B}))$.
Then $P^{B}_{t} f, \;  R^{B}_{\alpha}f$ is an $m$-version of $T^{B}_t f,  G^{B}_{\alpha}f$, respectively for any $f \in L^2(B,m)_b$. Since $P_t^{B} 1_{A}(x) \le P_t 1_{A}(x)$ for any $ A\in \mathcal{B}(B)$, $x \in B$ and $m$ has full support on $E$, $A \mapsto P_t^{B} 1_{A}(x), \; A \in \mathcal{B}(B)$ is absolutely continuous with respect to $1_{B} \cdot m$. Hence there exists a (measurable) transition kernel density $p^{B}_{t}(x, y)$, $x,y \in B$, such that
\begin{equation}\label{abspart}
P_t^{B} f(x) = \int_{B} p_t^{B} (x,y) \,  f(y) \, m(dy) ,\; t>0 \;, \;\; x\in B
\end{equation}
for $f  \in  \mathcal{B}_{b}(B)$. Correspondingly, there exists a (measurable) resolvent kernel density $r_{\alpha}^{B}(x,y)$, such that
\[
R_{\alpha}^{B} f(x) = \int_{B} r_{\alpha}^{B} (x,y) \,  f(y) \, m(dy) \,, \;\; \alpha>0, \;\;x \in B
\]
for $f  \in \mathcal{B}_{b}(B)$.
For a signed Radon measure $\mu$ on $B$, let us define
\[
R_{\alpha}^{B} \mu (x) = \int_{B} r_{\alpha}^{B} (x,y) \, \mu(dy) \, , \;\; \alpha>0, \;\;x \in B
\]
whenever this makes sense.
The process defined by
\begin{equation}\label{PP}
X^{B}_t(\omega)=
\begin{cases}
X_t(\omega), \;\;\;\; 0\le t < D_{B^c} (\omega) \\
\Delta, \;\;\;\; t \ge D_{B^c} (\omega)
\end{cases}
\end{equation}
is called the part process corresponding to $\E^{B}$ and is denoted by $\bM|_{B}$.  $\bM|_{B}$ is a Hunt process on $B$ (see \cite[p.174 and Theorem A.2.10]{FOT}).  In particular, by (\ref{abspart})
$\bM|_{B}$ satisfies the absolute continuity condition on $B$.\\ 
A positive Radon measure $\mu$ on $B$ is said to be of finite energy integral if
\[
\int_{B} |f(x)|\, \mu (dx) \leq C \sqrt{\E^{B}_1(f,f)}, \; f\in D(\E^{B}) \cap C_0(B),
\]
where $C$ is some constant independent of $f$. A positive Radon measure $\mu$ on $B$ is of finite energy integral (on $B$) if and only if there exists a unique function $U_{1}^{B} \, \mu\in D(\E^{B} )$ such that
\[
\E^{B}_{1}(U_{1}^{B} \, \mu, f) = \int_{B} f(x) \, \mu(dx),
\]
for all $f \in D(\E^{B}) \cap C_0(B)$. $U_{1}^{B} \, \mu$ is called $1$-potential of $\mu$. In particular, $R_{1}^{B} \mu$ is a version of $U_{1}^{B} \mu$ (see e.g. \cite[Exercise 4.2.2]{FOT}). The measures of finite energy integral are denoted by $S_0^{B}$. We further define $S_{00}^{B} : = \{\mu\in S_0^{B} \, | \; \mu(B)<\infty, \|U_{1}^{B} \mu\|_{\infty, B}<\infty \}$.
A positive Borel measure $\mu$ on $B$ is said to be smooth in the strict sense if there exists a sequence $(E_k)_{k\ge 1}$ of Borel sets increasing to $B$ such that $1_{E_{k}} \cdot \mu \in S_{00}^{B}$ for each $k$ and
\[
\P_{x} ( \lim_{k \rightarrow \infty} \sigma_{ B \setminus E_{k} }  \ge \zeta ) =1 \;, \;\; \forall x \in B.
\]
The totality of the smooth measures in the strict sense is denoted by $S_{1}^{B}$ (see \cite{FOT}). If  $\mu \in S_{1}^{B}$,
then there exists a unique $A \in A_{c,1}^{+, B}$ with $\mu = \mu_{A}$, i.e. $\mu$ is the Revuz measure of $A$ (see \cite[Theorem 5.1.7]{FOT}), such that
\[
\EE_x \Big[ \int_{0}^{\infty} e^{-t} \, d A_{t}  \Big] = R_1 \mu_{A}(x) \, , \;\; \forall x \in B.
\]
Here, $A_{c,1}^{+,B}$ denotes the positive continuous additive functionals on $B$ in the strict sense. If $B = E$, we omit the superscript $B$ and simply write $U_{1}, S_0, S_{00}, S_{1}$, and $A_{c,1}^{+}$.

\begin{lemma}\label{l;sumadd}
For $k \in \mathbb{Z}$, let $ \mu_{A^{k}}, \mu_{A} \in S_{1}^{B}$ be the Revuz measures associated with $A^{k}, A \in A_{c,1}^{+, B}$, respectively.
Suppose that $ \mu_{A} = \sum_{k \in \mathbb{Z}}  \mu_{A^{k}}$. Then $A= \sum_{k \in \mathbb{Z}} A^{k}$.
\end{lemma}
\pf
Since $\mu_{\sum_{-n \le k \le n } A^{k}} \le \mu_{A}$ and $\sum_{-n \le k \le n } A^{k} \in A_{c,1}^{+, B}$, we can use \cite[IV. (2.12) Proposition]{BlGe} in order to show that for any $n \in \N$ and $t \ge 0$
$$
\P_x \Big( \sum_{-n \le k \le n } A^{k}_{t} \le A_{t} \Big) = 1
$$
for all $x \in B$. Thus by the Weierstrass M-test
$\tilde{N_{t}} : =  \sum_{k \in \mathbb{Z}} A^{k}_{t}$ converges locally uniformly $\P_{x}$-a.s. for all $x \in B$. It follows that $\tilde{N_{t}}$ is positive continuous additive functional in the strict sense. In particular $d \tilde{N_{t}} = \sum_{k \in \mathbb{Z}} d A_t^k$ which further implies that for any $x \in B$ and $f \in C_{0}(B)$
\begin{eqnarray*}
&&\EE_x \Big[ \int_0^{\infty} e^{-t} \, f(X_t) \, d \tilde{N_{t}}  \Big] = \sum_{k \in \mathbb{Z}} \, \EE_x \Big[ \int_0^{\infty} e^{-t} \, f(X_t) \, d A_{t}^k  \Big]
= \sum_{k \in \mathbb{Z}} R_{1}^{B} f \mu_{A^k}(x) \\
&=& \sum_{k \in \mathbb{Z}} \int_{E} r_{1}^{B}(x,y)\, f(y) \; \mu_{A^{k}}(dy) =    \int_{E} r_{1}^{B} (x,y)\, f(y) \; \mu_{A}(dy) = R_{1}^{B} f  \mu_{A}(x) \\
&=& \EE_x \Big[ \int_0^{\infty} e^{-t} \, f(X_t) \, d A_{t}  \Big].
\end{eqnarray*}
Hence, $\tilde{N} = A$ by \cite[IV. (2.12) Proposition]{BlGe}.
\qed

\begin{proposition}\label{p;smooth}
Let $\mu$ be a positive Radon measure on $E$. Suppose that for some relatively compact open set $G \subset E$, $1_{G} \cdot \mu \in S_{0}$ and that $R_1(1_{G} \cdot \mu)$ is bounded $m$-a.e. on $E$ by a continuous function $r_1^G \in C(E)$ (resp. that $R_1 (1_{G} \cdot \mu) \in L^1(G,\mu)$ and that $R_1(1_G \cdot \mu)$ is bounded $m$-a.e. on $E$ by a continuous function $r_1^G \in C(E)$). Then $1_{G} \cdot \mu \in S_{00}$. In particular, if this holds for any relatively compact open set $G$, then $\mu \in S_{1}$ with respect to a sequence of open sets $ (E_{k})_{k \ge 1}$.
\end{proposition}
\pf
First suppose $1_G \cdot \mu \in S_0$. Since $\mu$ is a Radon measure, we have that $1_G \cdot \mu$ is finite. Since $r_1^G$ is continuous, it follows that
\[
E_k : = \{ x \in E \ | \  r_1^G(x) < k \}, \quad k \ge 1
\]
are open sets increasing to $E$. Let $\tilde{U}_1(1_{E_k \cap G} \cdot \mu)$, $\tilde{U}_1(1_{G} \cdot \mu)$ be q.c. versions of $U_1(1_{E_k \cap G} \cdot \mu)$, $U_1(1_{G} \cdot \mu)$. On $E_k$ it holds $\tilde{U}_1(1_{E_k \cap G} \cdot \mu) \le \tilde{U}_1(1_G \cdot \mu) \le r_1^G \le k$ q.e. 
Hence $U_1 ( 1_{E_k \cap G} \cdot \mu) \le k$ $m$-a.e. by \cite[Lemma 2.2.4 (ii)]{FOT}. Since $(E_k)_{k \ge 1}$ is an open cover of $\overline{G}$, we know that there exists $k_{0} \in \N$ with $G \subset \overline{G} \subset E_{k_{0}}$. Hence, 
 $U_1 \left(1_{G} \cdot \mu \right) \le k_{0}$ $m$-a.e. Therefore, $1_{G} \cdot \mu  \in S_{00}$. If $R_1(1_G \cdot \mu) \in L^1(G,\mu)$, then
\[
\int_G \int_G r_1(x,y) \ \mu(dy) \ \mu(dx) = \int_G  R_1(1_G \cdot \mu)(x) \ \mu(dx)< \infty. 
\]
Hence $1_G \cdot \mu \in S_0$ by \cite[Example 4.2.2]{FOT} and we conclude as before.
\qed

\section{Muckenhoupt weights}\label{3}

In this section we complete and extend substantially the results from \cite{ShTr13b}. We assume throughout that $E = \R^d$, with $d \ge 3$ 
(except in Lemma \ref{l;Feller}(vi),  Proposition \ref{p;3.8}(ii), Theorem \ref{t;3.9}(ii) and Remark \ref{r;3.15} where the state space is $\R^d\setminus \{0\}$ with $d\ge 2$). 
We consider a weight function that is in the Muckenhoupt $A_2$ class. For the definition and basic properties of Muckenhoupt weights, we refer to \cite{Tu}. Precisely, we assume the following:\\

($\alpha$) $\phi :\R^{d} \rightarrow [0,\infty)$  is a $\mathcal{B}(\R^d)$-measurable function and $\phi > 0$ $dx$-a.e.,\\

($\beta$) $\rho \phi \in  A_{2}, \;\;  \rho \in H^{1,1}_{loc}(\R^d, dx)$, $\rho > 0$ $dx$-a.e.,\\ 

and consider 
\begin{equation}\label{DF}
\E (f,g) : = \frac{1}{2}\int_{\R^d}\nabla f\cdot\nabla g \   dm, \quad f,g \in C_0^{\infty}(\R^d), \quad m := \rho \phi dx
\end{equation}
in $L^2(\R^d, m)$.

\begin{remark}\label{r;re1}
Let $\tilde{c} \ge 1$. If $\phi$ is measurable with $\tilde{c}^{-1} \le \phi \le \tilde{c}$ and $\rho \in A_{2}$, then $\rho \phi \in A_{2}$.
\end{remark}

Since $\rho \phi \in A_2$, we have $\frac{1}{\rho \phi}  \in L^1_{loc}(\R^d, dx)$, and the latter implies that  $\eqref{DF}$  is closable in $L^2(\R^d, m)$ (see \cite[II.2 a)]{MR}). The closure $(\E,D(\E))$ of $\eqref{DF}$ is a strongly local, regular, symmetric Dirichlet form (cf. e.g. \cite[p. 274]{St3}).\\
From \cite[p. 303 Proposition 2.3]{St2} and \cite[p. 286 A)]{St3} (see also \cite[5.B]{St3} and \cite{BM}) we know that there exists a jointly continuous transition kernel density  $p_{t}(x,y)$ such that
\[
P_t f(x) := \int_{\R^d} p_t(x,y)\, f(y) \ m(dy), \ \  t>0, \ x,y \in \R^d, \ f\in \mathcal{B}_b(\R^d)
\]
is an $m$-version of $T_t f$ if $f  \in  L^2(\R^d , m)_b$. We want to show that $(P_t)_{t \ge 0}$ is strong Feller. For this, we first need a lemma.

\begin{lemma}\label{l;majorant}
Let $t,r > 0$. Then $\inf_{x \in \bar{B}_r} m\big(B_{\sqrt{t}}(x)\big) = : M_{t,r} >0$ and for any $x \in \bar{B}_r$, $\varepsilon>0$
\begin{equation}\label{majorant}
p_t(x,y) \le  \frac{c \exp\Big(-\frac{\|y\|^2}{2(4+\varepsilon) t} \Big) \ \left( 1+\frac{\|y\|}{\sqrt{t}} \right)^{\alpha/2} }{M_{t,r}^{1/2} \ m \big(B_{\sqrt{t}}(0) \big)^{1/2}}  \ 1_{\R^d \setminus \bar{B}_{4r}} (y) + \Big( \sup_{\begin{subarray} \ x \in \bar{B}_r \\ y \in \bar{B}_{4r} \end{subarray}} p_t(x,y) \Big) 1_{\bar{B}_{4r}}(y)
\end{equation}
where $\alpha > 0$ is some constant. In particular
\[
\sup_{x \in \bar{B}_r} p_t(x, \cdot) \in L^1 (\R^d,m).
\]
\end{lemma}
\pf
It follows from \cite[4.3]{St4} and \cite[Corollary 4.2.]{St3} that  for $x, y \in \R^d$, $t > 0$ and any $\varepsilon >0$
\begin{equation}\label{heatesti}
p_t(x,y) \le c \ \frac{\exp{ \Big( -\frac{\|x - y\|^2}{(4+\varepsilon)t} \Big)}}{m \big(B_{\sqrt{t}} (x) \big)^{1/2} m \big(B_{\sqrt{t}} (y) \big)^{1/2}}.
\end{equation}
By Fatou's lemma, $x \mapsto m \big(B_{\sqrt{t}}(x) \big)$ is lower semicontinuous and so it attains its infimum on $\bar{B}_r$. Therefore $M_{t,r}>0$. Moreover, since $- \| x-y\|^2 \le - \frac{\|y\|^2}{2} + \frac{\| y\| (4 \|x\| - \| y\|)}{2}$, we obtain $- \| x-y\|^2 \le - \frac{\|y\|^2}{2} $ for any $x \in \bar{B}_r$ if $y \in \R^d \setminus \bar{B}_{4r}$. Further for some $\alpha >0$ and any $x,y \in \R^d$
\begin{equation}\label{15a}
m \big(B_{\sqrt{t}}(y)\big) \ge \frac{m(B_{\sqrt{t}}(x))  }{C_D} \left(1 + \frac{\|x-y\|}{\sqrt{t}} \right)^{-\alpha},
\end{equation}
where $C_D$ is the volume doubling constant of $m$ (see \cite[Proposition 5.1]{GrHu}). These facts together with the joint continuity of $p_t(x,y)$ and \eqref{heatesti} lead to \eqref{majorant}. Since $m \big(B_{r}(y) \big)$ has at most polynomial growth in $r$ for any $y \in \R^d$ (cf. Proof of Proposition 2.4 in \cite{ShTr13b}) the last statement follows. 
\qed

\begin{proposition}\label{p;strongfl}
\begin{itemize}
\item[(i)] $(P_t)_{t \ge 0}$ (resp. $(R_{\alpha})_{\alpha > 0}$) is strong Feller, i.e. for $t>0$, we have $P_t(\mathcal{B}_b(\R^d)) \subset C_b(\R^d)$ (resp. for $\alpha>0$, we have $R_{\alpha}(\mathcal{B}_b(\R^d)) \subset C_b(\R^d)$).
\item[(ii)] $(\textbf{H1})$ and \textbf{(H2)}$^{\prime}$ (iii) and (iv) hold for $(P_t)_{t \ge 0}$.
\item[(iii)] $P_t (L^1(\R^d,m)_0) \subset C_{\infty} (\R^d)$.
\item[(iv)] Let $\mu$ be a positive Radon measure and $G \subset \R^d$ relatively compact open. Let
\[
\int_G r_1(\cdot,y) \ \mu(dy) \le r_1^G
\]
$\mu$-a.e. on $G$ and $m$-a.e. on $\R^d$, where $r_1^G$ is a continuous function on $\R^d$. Then $1_G \cdot \mu \in S_{00}$.
\end{itemize}
\end{proposition}
\pf
(i) Let $x_n \to x$ in $\R^d$. For  $f \in \mathcal{B}_b(\R^d)$ and $t > 0$
\[
| P_t f(x_n) - P_t f(x) | \le \int_{\R^d} | p_t(x_n,y) - p_t(x,y) | \ | f(y) | \ m(dy)
\]   
which converges to $0$ by Lebesgue in view of Lemma \ref{l;majorant} and the continuity of $p_t(\cdot,y)$. Clearly, $P_t f $ is bounded. Hence, $(P_t)_{t \ge 0}$ is strong Feller. Since $R_{\alpha} f (x) = \int_0^{\infty} e^{-t} \ P_t f(x) \ dt$ and $\|P_t f\|_{\infty} \le \|f\|_{\infty}$ for any $f \in \mathcal{B}_b(\R^d)$, $(R_{\alpha})_{\alpha>0}$ is clearly also strong Feller by Lebesgue.\\
(ii) By (i), $A \mapsto P_t(x,A)$ is a sub-probability measure on $\mathcal{B}(\R^d)$ for any $t>0$, $x \in \R^d$. Obviously, $x  \mapsto P_t(x,A)$ is also measurable for any $A \in \mathcal{B}(\R^d)$ and so it remains to show the Chapman-Kolmogorov equation.
 By the semigroup property, 
\begin{equation}\label{semiproper}
P_{t+s} 1_A(x) = P_t ( P_s 1_A ) (x), \ \ A \in \mathcal{B}(\R^d), \ t,s >0
\end{equation}
for $m$-a.e.$x \in \R^d$. From the strong Feller property, both sides of \eqref{semiproper} are continuous, hence \eqref{semiproper} holds for every $x \in \R^d$, i.e. the Chapman-Kolmogorov equation holds and so $(P_t)_{t \ge 0}$ is a  sub-Markovian transition function. $\textbf{(H2)}^{\prime}$ (iii) follows from (i) and $\textbf{(H2)}^{\prime}$ (iv) follows from \cite[Proposition 3.1]{St3}.\\
(iii) Combining \eqref{heatesti} and \eqref{15a} we have for any $x,y \in \R^d$, $t > 0$ and $\varepsilon > 0$,  
\begin{equation}\label{densest2}
p_t(x,y) \le c \frac{1}{m \left(B_{\sqrt{t}} (y) \right)} \ \exp{ \Big( -\frac{\|x - y\|^2}{(4+\varepsilon) t} \Big)}. 
\end{equation}
Using the joint continuity of $p_t(\cdot,\cdot)$, as in (i) we can see that $P_t(L^1(\R^d,m)_0) \subset C(\R^d)$. Let $f \in L^1(\R^d,m)_0$. Using \eqref{densest2}, 
\[
\big| P_t f(x) \big| 
\le \frac{c}{\inf_{y \in \text{supp}(f)} m(B_{\sqrt{t}}(y))} \int_{\text{supp}(f)} |f(y)| \ e^{-\frac{\|x-y\|^2}{(4+ \varepsilon)t}} \ m(dy)
\]
which converges to 0 by Lebesgue as $x \to \infty$.\\
(iv) This is just a reformulation of Proposition \ref{p;smooth}. 
\qed

First let us assume that\\

($\gamma$) The transition function $(P_t)_{t \ge 0}$ satisfies (\textbf{H2}) with $E= \R^d$. \\

Later we will use the Feller semigroup method and the Dirichlet form method for some typical Muckenhoupt $A_2$ weights to verify ($\gamma$). By the existence of $\bM$ associated with $(P_t)_{t \ge 0}$, $\bM$ satisfies the absolute continuity condition. Since $\rho \phi \in A_2$,
$(\E,D(\E))$ is conservative, i.e. $T_t 1(x)=1$ for $m$-a.e. $x\in \mathbb{R}^d$ and all $t>0$ (see \cite[Proposition 2.4]{ShTr13b}). It follows
\begin{equation}\label{HPcons}
\P_{x}(\zeta = \infty)=1, \ \ \forall x \in \R^d,
\end{equation}
by \cite[Theorem 4.5.4 (iv)]{FOT} and 
\begin{equation}\label{eq20}
\P_x \big(t \mapsto X_{t} \text{ is continuous on} \  [0, \infty)\big) =1, \ \ \forall x \in \R^d,
\end{equation}
by \cite[Theorem 4.5.4 (ii)]{FOT}.\\
Throughout, let  $f^j (x):= x_j$, $j=1,\dots, d$, $x \in \R^d$, be the coordinate projections. In order to be explicit, we further assume the following integrations by parts formula\\

(IBP) For $ f  \in \{f^1, \dots,f^d\}$, $g \in C_{0}^{\infty}(\R^{d})$\\

\[
-\E (f , g) = \int_{\R^d} \left (\nabla f  \cdot  \frac{\nabla \rho  }{2\rho}  \right ) g \ dm
+ \int_{\R^{d}} g \, d\nu^{f},
\]
where $\nu^{f}= \sum_{k \in \mathbb{Z}} \nu_{k}^{f}$ and $\nu^{f},  \nu_{k}^{f}, k \in \mathbb{Z}$ are signed Radon measures (locally of bounded total variation).\\

For a signed Radon measure $\mu$ we denote by $\mu^+$ and $\mu^-$ the positive and negative parts in the Hahn decomposition for $\mu$, i.e. $\mu = \mu^+ -\mu^-$.
Additionally, we assume that\\

($\delta$) For any $G \subset \R^d$ relatively compact open, $k \in \mathbb{Z}$ and $f \in \{f^1, \dots, f^d\}$, we have that $1_G \cdot \nu^{f+}$, $1_G \cdot  \nu^{f-}$, $1_G  \cdot  \nu_k^{f+}$, $1_G \cdot  \nu_k^{f-}, 1_G \cdot  \frac{\|\nabla \rho\|}{\rho} m \in S_0$ and the corresponding 1-potentials are all bounded by continuous functions.\\

\begin{thm}\label{t;Ldec}
Suppose $(\alpha)-(\delta)$ and (IBP). Then
\begin{equation}\label{LSFD1}
X_{t} = x + W_{t} + \int_{0}^{t}\frac{\nabla \rho}{2 \,\rho}(X_s) \, ds + \sum_{k \in \mathbb{Z} } L_{t}^{k} \;,\;\; t  \ge 0,
\end{equation}
$\P_x$-a.s. for any $x \in \R^{d}$ where $W$ is a standard d-dimensional Brownian motion starting from zero, $L^{k} = (L^{1, k}, \dots, L^{d, k})$ and $L^{j,k}, j=1,\dots, d$, is the difference of positive continuous additive functionals of $X$  in the strict sense associated with Revuz measure $\nu_{k}^{f^{j}} = \nu_{k}^{f^{j},(1)} - \nu_{k}^{f^{j},(2)}$ defined in (IBP) (cf. \cite[Theorem 5.1.3]{FOT}).
\end{thm}
\pf
Given that $(\alpha)-(\delta)$ and (IBP) hold, the assertion follows from \cite[Theorems 5.1.3 and 5.5.5]{FOT}, Lemma \ref{l;sumadd}, and Propositions \ref{p;strongfl} and \ref{p;smooth}.
\qed

For later purpose we add some auxiliary results. Define 
\begin{equation}\label{rieszpot}
V_{\eta} g(x) : = \int_{\R^d} \frac{1}{\| x-y \|^{d-\eta}} \ g(y) \ dy, \quad x \in \R^d, \ \eta > 0,
\end{equation}
whenever it makes sense.
\begin{lemma}\label{l;miz}
Let $\eta \in (0,d)$, $0<\eta - \frac{d}{p} <1$ and $g \in L^p(\R^d, dx)$ with  
\[
\int_{\R^d} (1+\|y\|)^{\eta -d} |g(y)| \ dy < \infty.
\]
Then $V_{\eta}g$ is H\"older continuous of order $\eta - \frac{d}{p}$.
\end{lemma}
\pf
See \cite[Chapter 4, Theorem 2.2]{Miz}.
\qed

\begin{lemma}\label{l;Feller}
Let $\tilde{c}^{-1} \|x\|^{\alpha} \le \rho \phi(x) \le \tilde{c} \|x\|^{\alpha}$ for some $\alpha \in (-d,d)$, $\tilde{c} \ge 1$.
Then:
\begin{itemize}
\item[(i)] $\lim_{t \downarrow 0}P_t f(x) = f(x)$, $\forall x \in \R^d$, $\forall f \in C_0(\R^d)$, i.e. (\textbf{H1}) and (\textbf{H2}) hold (cf. Proposition \ref{p;strongfl}(i),(iii) and Lemma \ref{t;Feller}).
\item[(ii)] Let $\Phi(x,y): = \frac{1}{\|x-y\|^{\alpha +d -2}}$ and $\Psi(x,y):= \frac{1}{\|x-y\|^{d -2} \|y\|^{\alpha}}$. Then 
\[
c^{-1} \left(\Phi(x,y) + \Psi(x,y)1_{\{ \alpha \in [0,d)\}} \right)  \le r_1(x,y) \le c \left( \Phi(x,y) + \Psi(x,y)1_{\{ \alpha \in (-d,0)\}}  \right).
\]
\item[(iii)] Let $\alpha \in (-d+1,2) $ and $G\subset \R^d$ any relatively compact open set. Suppose $1_G  \cdot f \ \| x\|^{\alpha} \in L^p(\R^d,dx)$, $p \ge 1$ with $0 < 2- \alpha - \frac{d}{p} <1$ and $1_G \cdot  f \in L^q(\R^d,dx)$ with $0 < 2- \frac{d}{q} <1$. Then $R_1(1_G \cdot  |f| m)$ is bounded everywhere (hence clearly also bounded $m$-a.e. on $\R^d$ and $R_1(1_G |f| m) \in L^1(G,|f|m)$) by the continuous function 
$\int_G |f(y)| \ \left(\Phi(\cdot,y) + \Psi(\cdot,y) \right) \ m(dy)$. In particular, Proposition \ref{p;smooth} applies and $1_G  \cdot  |f| m \in S_{00}$.
\item[(iv)] Let $\alpha \in (-d+1,2)$. Then $R_1\left(1_G \cdot  \frac{\|\nabla \rho\|}{\rho} m \right)$ is pointwise bounded by a continuous function for any relatively compact open set $G \subset \R^d$. In particular $1_G \cdot  \frac{\|\nabla \rho\|}{\rho} m \in S_{00}$ for any relatively compact open set $G \subset \R^d$.
\item[(v)] Let  $\alpha \in (-d+1,1)$. Let $D \subset \R^d$ be a bounded Lipschitz domain with surface measure $\sigma_{\partial D}$. Suppose that $\rho$ is bounded on $\partial D$ (more precisely the trace of $\rho$ on $\partial D$, which exists since $\rho \in H^{1,1}_{loc}(\R^d)$). Then $R_1(1_G \cdot \rho  \sigma_{\partial D})$ is pointwise bounded by a continuous function for any relatively compact open set $G \subset \R^d$. In particular $1_{G} \cdot \rho  \sigma_{\partial D} \in S_{00}$ for any relatively compact open $G \subset \R^d$.
\item[(vi)] Let $\alpha \in [-d+2,d)$, $d\ge 2$. Then $\emph{Cap}(\{0\})=0$ and the part Dirichlet form $(\E^B,D(\E^B))$ on $B:=\R^d \setminus \{0\}$ satisfies (\textbf{H1}), (\textbf{H2}) with transition kernel density $p_t^B = p_t |_{B \times B}$. Moreover $(\E^B,D(\E^B))$ is conservative.
\end{itemize}
\end{lemma}
\pf
(i) From Proposition \ref{p;strongfl}, we know that $(P_t)_{t \ge 0}$ satisfies (\textbf{H1}) and is strong Feller. Thus Lemma \ref{t;Feller} (ii) holds. We will check Lemma \ref{t;Feller} (i).
Let $m_{\alpha} : = \|y\|^{\alpha}dy$, $\alpha \in (-d,d)$. For $\alpha \in [0,d)$ and $0< \sqrt{t} \le \|x\|$, we have 
\begin{eqnarray}\label{20}
m(B_{\sqrt{t}}(x)) \ge c_d \tilde{c}^{-1} \ (\| x\| - \sqrt{t})^{\alpha} \ \sqrt{t}^d,
\end{eqnarray}
with $c_d = vol(B_1(0))$, and for $\alpha \in (-d,0)$ and $0 < \sqrt{t} \le \| x\|$, we have
\begin{equation}\label{21}
m(B_{\sqrt{t}}(x)) \ge \tilde{c}^{-1} c_d \ \sqrt{t}^{d} \ (2 \| x \|)^{\alpha}.
\end{equation}

Since $(\E,D(\E))$ is conservative and $(P_t)_{t \ge 0}$ is strong Feller, we have $P_t 1(x) = 1$ for all $x \in \R^d$, $t >0$. Thus by \eqref{20}, symmetry of $p_t(x,y)$ in $(x,y)$, and \eqref{densest2}, we get
\[
\big| P_t f(x) - f(x) \big| \le c \int_{\R^d} \big| f(x+\sqrt{t} y) - f(x) \big| \,\frac{ \exp \left( - \frac{\|y\|^2}{4+ \varepsilon} \right) \   }{  (\|x\| - \sqrt{t})^{\alpha}   } \ \|x + \sqrt{t} y\|^{\alpha}   \  dy,
\]
which converges to $0$ as $t \to 0$. For $x=0$, by \eqref{densest2} and symmetry of $p_t(\cdot,\cdot)$, we get 
\[
\big| P_t f(0) - f(0) \big| \le  c \int_{\R^d} \big| f(\sqrt{t} y) - f(0) \big| \  \exp \left( - \frac{\|y\|^2}{4+ \varepsilon} \right) \,  \ \|  y\|^{\alpha}  \ dy,
\]
which also converges to $0$ as $t \to 0$.
For $\alpha \in (-d,0)$, using \eqref{21} instead of \eqref{20}, similarly to the case of $\alpha \in [0,d)$ one can show that $P_t f(x) \to f(x)$ as $t \to 0$. Thus Lemma \ref{t;Feller} (i) holds.\\
(ii)
For $\alpha \in [0,d)$, we have
\begin{equation}\label{vol21}
 c  \ \sqrt{t}^{\alpha + d} \le m(B_{\sqrt{t}}(x)) \le c  \ \sqrt{t}^{d} \ (\| x \| + \sqrt{t})^{\alpha},
\end{equation}
and for $\alpha \in (-d,0)$, 
\begin{equation}\label{vol22}
c \ \sqrt{t}^{d} \ (\| x \| + \sqrt{t})^{\alpha} \le m \left(B_{\sqrt{t}} (x) \right) \le c \sqrt{t}^{d + \alpha}.
\end{equation}
By \cite[Corollary 4.10]{St3} and \eqref{densest2} 
\[
\frac{1}{c \ m \left(B_{\sqrt{t}} (y) \right)} \ \exp{ \Big( -c \frac{\|x - y\|^2}{t} \Big)} \le p_t(x,y) \le  \frac{c}{m \left(B_{\sqrt{t}} (y) \right)} \ \exp{ \Big( - \frac{\|x - y\|^2}{5t} \Big)}.
\]
Let first $\alpha \in [0,d)$. Then, for $x,y \in \R^d$ using the first inequality in \eqref{vol21}, we get
\[
r_1(x,y)  \le  \int_0^{\infty}  \ \frac{c}{(\sqrt{t})^{\alpha+d}} \exp \Big( - \frac{\|x-y\|^2}{5t} \Big)  \ dt.
\]
By standard calculations, using a change of variable with $s= \frac{\|x-y\|^2}{t}$, we obtain
\begin{equation}\label{resest1}
r_{1} (x,y)  \leq  \frac{c}{\|x-y\|^{\alpha + d-2}}.
\end{equation}
Using the second inequality in \eqref{vol21}, we get the lower bound of $r_1(x,y)$,
\begin{eqnarray*}
r_1(x,y) &\ge& \int_0^{\|y\|^2}  \ \frac{c}{\sqrt{t}^{d} \ ( 2 \|y\|)^{\alpha}} \exp \Big( - c \frac{\|x-y\|^2}{t} \Big)  \ dt \\
&& + \int_{\|y\|^2}^{\infty}  \ \frac{c}{\sqrt{t}^{d} \ (2 \sqrt{t} )^{\alpha}} \exp \Big( - c \frac{\|x-y\|^2}{t} \Big)  \ dt.
\end{eqnarray*}
Hence,
\[
r_{1} (x,y)  \ge  c \left( \frac{1}{\|x-y\|^{\alpha +d -2}} + \frac{1}{\|x-y\|^{d -2} \|y\|^{\alpha}} \right).
\]
For $\alpha \in (-d,0)$, using \eqref{vol22} instead of \eqref{vol21}, we get
\[
\frac{c}{\|x-y\|^{\alpha + d-2}} \le r_{1} (x,y)  \le  c \left( \frac{1}{\|x-y\|^{\alpha +d -2}} + \frac{1}{\|x-y\|^{d -2} \|y\|^{\alpha}} \right).
\]
(iii)
If $\alpha \in (-d+2,2)$, then the conditions imply that $V_{2-\alpha} \left(1_G \cdot  f \| x \|^{\alpha}\right)$ as well as $V_{2}\left(1_G \cdot f\right)$ are continuous by Lemma \ref{l;miz}. If $\alpha \in (-d+1,-d+2]$, then $V_{2-\alpha}\left(1_G \cdot f \| x \|^{\alpha}\right)$ is easily seen to be also continuous and so by (ii) for any $x \in \R^d$
\[
R_1(1_G \cdot |f|m)(x) \le c \ \big( \ V_{2-\alpha} (1_G \cdot |f| \ \|x\|^{\alpha})(x) + V_2(1_G \cdot |f|)(x) \big).
\]
(iv) 
Let $\alpha \in (-d+2,2)$ and $0 < \varepsilon < 1$ satisfy $2-\varepsilon > \alpha$. Then $1_G \cdot \frac{\|\nabla \rho\|}{\rho} \|x\|^{\alpha} = c \ 1_G \cdot \|x\|^{\alpha -1} \in L^p(\R^d,dx)$ for $p = \frac{d}{(2-\varepsilon)-\alpha} \ge 1$ and $1_G \cdot \frac{\|\nabla \rho\|}{\rho} = c \ 1_G \cdot \|x\|^{-1} \in L^{q}(\R^d,dx)$ for $q= \frac{2d}{3}$.  For $\alpha \in (-d+1,-d+2]$, $V_{2-\alpha}\left(1_G \cdot \frac{\|\nabla \rho\|}{\rho} \|x\|^{\alpha}\right)$ is continuous. Therefore, the statement follows as in (iii).\\
(v) Let $G$ be relatively compact open. There exist $B_i \subset \partial D$ and Lipschitz continuous functions $F_i$, $i=1, \dots n$ for some $n \in \N$ such that $B_i = \{ x \in \partial D \ | \ x = (x',F_i(x')) \in \R^{d-1} \times \R \}$ and $\bigcup_{i=1}^{n} B_i = \partial D$. Set $B_i^* = \{y' \in \R^{d-1} \ | \ (y',F_i(y')) \in B_i )  \}$. Let first $\alpha \in [0,1)$. Then, using (ii) we get for every $x \in \R^d$ 
\begin{eqnarray*}
R_1 (1_G \cdot \rho \sigma_{\partial D})(x) &\le& c \sum_{i=1}^{n} \int_{B_i^*} \frac{\sqrt{1+|\grad F_i(y')|^2}}{\|x'-y'\|^{\alpha + d-2}}  \ dy'\\
&\le& c \sum_{i=1}^{n} \int_{B_i^*} \frac{1}{\|x'-y'\|^{\alpha + d-2}}  \ dy'\\
&\le& c \int_K \frac{1}{\|x' - y'\|^{\alpha + d -2}} \ dy',
\end{eqnarray*}
where $K \subset \R^{d-1}$ is some compact set. Since the last expression is continuous in $x'$ (hence in particular in $x$) by Lemma \ref{l;miz}, the final statement follows by Proposition \ref{p;smooth}.
For $\alpha \in (-d+1,0)$ we have for any $x \in \R^d$ that 
\begin{eqnarray*}
R_1 (1_G \cdot \rho \sigma_{\partial D})(x) &\le& c \int_{\partial D} \left( \frac{1}{\|y\|^{\alpha} \|x-y \|^{d-2}} + \frac{1}{\|x-y \|^{\alpha + d -2}} \right) \ \| y\|^{\alpha} \ \sigma_{\partial D} (dy)\\ 
&\le& c\int_{\partial D}  \left( \frac{1}{ \|x-y \|^{d-2}} + \frac{1}{\|x-y \|^{\alpha + d -2}} \right)  \sigma_{\partial D} (dx),\\
\end{eqnarray*}
and we conclude as before in the case of $\alpha \in [0,1)$.\\
(vi) By \cite[Example 3.3.2]{FOT} it holds $\mbox{Cap}(\{0\})=0$. Hence $u(x):=\P_x(\sigma_{B^c} < \infty) =0$ for $m$-a.e. $x \in \R^d$. Since $u$ is an excessive function and $\bM$ satisfies the absolute continuity condition it follows $u(x)=0$ for all $x \in \R^d$. From this the remaining part of the proof is straightforward. 
\qed

\subsection{Skew reflection on spheres and on a Lipschitz domain}\label{3.1}
\subsubsection{Skew reflection on spheres}\label{3.1.1}
In \cite{ShTr13b}, we considered the Dirichlet form determined by \eqref{DF} with concrete $\phi$ and $\rho = \xi^2$, $\xi \in H^{1,2}_{loc}(\R^d)$. More precisely, our assumptions were the followings:
we let $m_0 \in (0, \infty)$ and $(l_{k})_{k \in \mathbb{Z}} \subset (0, m_0)$, $0 < l_{k} < l_{k+1} <m_0$, be a sequence converging to $0$ as $k \rightarrow -\infty$ and converging to $m_0$ as $k \rightarrow \infty$, $(r_{k})_{k \in \mathbb{Z}} \subset (m_0, \infty)$, $m_0 < r_{k} < r_{k+1} < \infty$, be a sequence converging to $m_0$ as $k \rightarrow -\infty$ and tending to infinity  as $k \rightarrow \infty$, and set
\begin{equation}\label{eqphi}
\phi : = \sum_{k \in \mathbb{Z} } \left (\gamma_{k}  \cdot 1_{A_{k}}
+ \overline{\gamma}_{k} \cdot 1_{\hat{A}_{k}}\right ),
\end{equation}
where $\gamma_{k}$ , $ \overline{\gamma}_{k} \in (0, \infty)$, $A_{k} : = B_{l_{k}} \setminus \overline{B}_{l_{k-1}}$ , $\hat{A}_{k} : = B_{r_{k}} \setminus \overline{B}_{r_{k-1}}$. 
We further assumed
\begin{itemize}
\item[(a)]  $\rho \phi \in A_{2}$, $\rho = \xi^2$, $\xi \in H^{1,2}_{loc}(\R^d)$, $\rho > 0$ $dx$-a.e.
\item[(b)]  $\sum_{k \in \mathbb{Z}} |\,  \gamma_{k+1} - \gamma_{k} \,| + \sum_{k \le 0} |\,  \overline{\gamma}_{k+1} - \overline{\gamma}_{k} \,| < \infty$ and for all $r>0$ there exists a $\delta_r > 0$ such that $\phi \ge \delta_r$ $dx$-a.e. on $B_r$.
\item[(c)] $(\gamma)$ is satisfied, $R_1 \big( L^2(\R^d,m)_0 \big) \subset C(\R^d)$, and if $\phi \nequiv 1$ then $R_1 (1_G \cdot \rho \sigma_r) \in C(\R^d)$ for any $G \subset \R^d$ relatively compact open and $r>0$, where $\sigma_r$ is the surface measure on $\partial B_r$.
\end{itemize}
Under the assumptions (a)-(c), we showed (see \cite[Theorem 2.6]{ShTr13b}) that $\bM$ satisfies
\begin{equation}\label{equast}
X_t=x+W_t+\int_0^t \frac{\nabla \rho}{2\rho}(X_s) \ ds+\int_0^{\infty}\int_0^t \nu_{a}(X_s) \, d\ell_s^a(\|X\|)\,\eta(da),\ t \ge 0,
\end{equation}
$\P_x$-a.s. for any $x \in \R^d$, where $W$ is a $d$-dimensional standard Brownian motion starting at $0$, $\nu_{a}$ 
is the unit outward normal on the boundary $\partial B_a$, $\ell^{a}(\|X\|)$ is the 
symmetric semimartingale local time at $a\in (0,\infty)$ of $\|X\|$, $\eta=\sum_{k\in \mathbb{Z}}(2\alpha_k-1)\delta_{d_k}$ with $(\alpha_k)_{k\in \mathbb{Z}}\subset (0,1)$ is a sum of Dirac measures at a sequence  $(d_k)_{k\in \mathbb{Z}}\subset (0,\infty)$ with exactly two accumulation points in $[0,\infty)$, one is zero and the other is $m_0>0$, and $(\alpha_k)_{k \in \mathbb{Z}}$ and $(d_k)_{k \in \mathbb{Z}}$ are determined by $(\gamma_k)_{k \in \mathbb{Z}}$, $(\overline{\gamma}_k)_{k \in \mathbb{Z}}$, $(l_k)_{k \in \mathbb{Z}}$, and $(r_k)_{k \in \mathbb{Z}}$ (see \cite{ShTr13b}).  
\begin{remark}\label{r;compar}
The assumptions (a)-(c) imply $(\alpha)-(\delta)$ and (IBP). However, in comparison to \cite{ShTr13b}, we insist to point out two improvements. The first one is that in $(\alpha)$ $\rho$ is only assumed to be in $H^{1,1}_{loc}(\R^d)$ instead of $\rho = \xi^2$ with  $\xi \in H^{1,2}_{loc}(\R^d)$ in (a). ($\alpha$) allows to consider weights that increase rapidly  toward singularities which are of positive capacity. A typical example is $\rho(x) = \|x\|^{\alpha}, \ \alpha \in (-d+1, -d+2]$ (cf. \cite[Example 3.3.2]{FOT}). The second improvement is that in ($\delta$) the potentials are only assumed to be bounded by continuous functions and not to be continuous as in $(c)$ (cf. Proposition \ref{p;smooth} and \cite[Lemma 4.2, Theorem 4.5]{ShTr13b}). In particular, replacing (a) with $(\beta)$, and (c) with $(\gamma)$ and $(\delta)$, we still obtain \eqref{equast}. 
\end{remark}

\begin{prop}\label{p;3.8}
\begin{itemize}
\item[(i)]Let $\rho(x) = \|x\|^{\alpha}$, $\alpha \in (-d+1,2)$. Let $\phi$ be like in \eqref{eqphi}, satisfy (b) and $\tilde{c}^{-1} \le \phi \le \tilde{c}$ for some $\tilde{c} \ge 1$. If $\phi \equiv \tilde{c}$ $dx$-a.e. (i.e. $\eta(da) \equiv 0$) or $\phi \nequiv \tilde{c}$ $dx$-a.e. and $\alpha \in (-d+1,1)$, then \eqref{equast} holds.
\item[(ii)] Let $\rho(x) = \|x\|^{\alpha}$, $\alpha \in [1,d)$, $d\ge 2$. Let $\phi$ be like in \eqref{eqphi}, satisfy (b) and $\tilde{c}^{-1} \le \phi \le \tilde{c}$ for some $\tilde{c} \ge 1$. Then \eqref{equast} holds for any $x \in \R^d \setminus \{0\}$. 
\end{itemize}
\end{prop}
\pf
(i) The proof follows from Proposition \ref{p;strongfl}, Lemma \ref{l;Feller} (i), (iv) and (v), and Remark \ref{r;compar}.\\
(ii) By Lemma \ref{l;Feller} (vi) $(\E^B,D(\E^B))$, $B:= \R^d \setminus \{0\}$ satisfies (\textbf{H1}), (\textbf{H2}). Fix $\alpha \in [1,d)$. Let 
\[
B_k : = \Big\{x \in \R^d \ \Big| \ \frac{l_{-k+1} + l_{-k}}{2} < \|x\| < \frac{r_{k+1} + r_k}{2} \Big\}, \quad k \ge 1.
\]
Then
\[
b_k := \tilde{c}^{-1} \Big(\frac{l_{-k+1} + l_{-k}}{2}\Big)^{\alpha} < \rho \phi  < \tilde{c} \ \Big(\frac{r_{k+1} + r_k}{2}\Big)^{\alpha} = : e_k
\]
on $B_k$. Set $d_k := \max (b_k^{-1}, e_k)$, $k \ge 1$. Then $(B_k)_{k \ge 1}$ is an increasing sequence of relatively compact open sets with smooth boundary such that $\bigcup_{k \ge 1} B_k = \bigcup_{k \ge 1} \overline{B}_k = B$ and $ \rho \phi \in (d_k^{-1},d_k)$ on $B_k$ where $d_k \to \infty$ as $k \to \infty$. Moreover $\|\nabla \rho \| \in L^{\infty}(B_k,dx)$ for any $k \ge 1$. We may now apply a localization procedure similarly to the one that is presented in all details subsequently to Lemma \ref{l;part3} to obtain the assertion. We only note that by the Nash type inequality of Lemma \ref{T;NI1} we obtain resolvent kernel estimates as in Corollary \ref{C;RDE1}. These local resolvent estimates are better than the global ones of Lemma \ref{l;Feller} (ii).
\qed

\subsubsection{Skew reflection on a Lipschitz domain}\label{3.1.2}

We consider the Dirichlet form determined by \eqref{DF} with $\rho(x) = \|x\|^{\alpha}$, $\alpha \in (-d+1,d)$ and
\begin{equation}\label{eqphi2}
\phi(x) : = \beta \, 1_{G^{c}}(x) + (1- \beta) 1_{G} (x)
\end{equation}
where $\beta \in (0,1)$ and $G \subset \R^d$ is a bounded Lipschitz domain. Then the following integration by parts formula holds for $f\in \{f^1,\dots,f^d\}$, $g \in C_0^{\infty}(\R^d)$
\[
-\E(f,g) = \int_{\R^d} \left(\nabla f \cdot \frac{\nabla \rho}{2 \rho} \right) \ g \ dm + (2\beta -1) \int_{\partial G} \nabla f \cdot \nu \ \frac{\rho } {2} \ d \sigma,
\]
where $\nu$ denotes the unit outward normal on $\partial G$ (cf. \cite{Tr1} and \cite{Tr3}). The existence of a Hunt process associated to $\E$ that satisfies the absolute continuity condition follows from Lemma \ref{l;Feller} (i). Furthermore:
\begin{thm}\label{t;3.9}
\begin{itemize}
\item[(i)] Let $\alpha \in (-d+1,1)$. Then
\begin{equation}\label{SFD1} 
X_{t} = x + W_{t} + \alpha \int_{0}^{t}\frac{X_s}{2 \|X_s\|^2} \, ds + (2 \beta -1) \int_0^t  \nu (X_s)\, d \ell_s \quad t  \ge 0
\end{equation}
$\P_x$-a.s. for all $x \in \R^d$, where $(W_t)_{t \ge 0}$ is a d-dimensional Brownian motion starting from zero and $(\ell_t)_{t\ge 0} \in A^{+}_{c,1}$ is uniquely associated to the surface measure $\frac{\rho \sigma}{2}$ on $\partial G$ via the Revuz correspondence.
\item[(ii)] Let $0 \notin \partial G$ and $\alpha \in [1,d)$, $d\ge 2$. Then \eqref{SFD1} holds $\P_x$-a.s. for any $x \in \R^d \setminus \{0\}$.
\end{itemize}
\end{thm}
\pf
(i)  Lemma \ref{l;Feller} (iv) and (v) apply. Therefore $(\alpha)$-$(\delta)$ and (IBP) are satisfied and the assertion immediately follows from Theorem \ref{t;Ldec}.\\
(ii) The proof is similar to the proof of Proposition \ref{p;3.8} (ii). We therefore only indicate the sequences $(B_k)_{k \ge 1}$ and $(d_k)_{k \ge 1}$. Fix $\alpha \in [1,d)$ We have either
$0 \in G$ or $0 \in \overline{G}^c$. If $0 \in G$, then choose $k_0 \ge 1$ such that $\partial G \subset \{x \in \R^d \ | \ k_0^{-1} < \|x\| < k_0\}$ and let
\[
B_k := \{x \in \R^d \ | \ (k_0 + k)^{-1} < \|x\| < k_0 + k \}, \quad  k \ge 1.
\]
Then
\[
b_k : = \min(\beta, 1-\beta)(k_0 + k)^{-\alpha} <  \rho \phi < \max(\beta,1-\beta)(k_0 + k)^{\alpha} = : e_k
\]
and we let $d_k : = \max(b_k^{-1},e_k)$, $k \ge 1$. If $0 \in \overline{G}^c$ then similarly we can find suitable $(B_k)_{k \ge 1}$ and $(d_k)_{k \ge 1}$.
\qed
\begin{remark}
This result was announced in \cite{ShTr13b} and extends a result obtained by Portenko in \cite[III, \S 3 and \S 4]{port}.
\end{remark}

\subsection{Skew reflection on hyperplanes}\label{sub3}

We consider skew reflection on hyperplanes
\[
H_{s} := \{ x \in \R^d \,|\ x_d =s \}, \;\; s \in \R.
\]
Let  $(l_{k})_{k \in \mathbb{Z}} \subset (-\infty, 0)$, $- \infty < l_{k} < l_{k+1} < 0$  be a sequence converging to $0$ as $k \rightarrow \infty$ and tending to $- \infty$ as $k \rightarrow -\infty$. Let $(r_{k})_{k \in \mathbb{Z}} \subset (0, \infty)$, $0 < r_{k} < r_{k+1} < \infty$ be a sequence converging to $0$ as $k \rightarrow -\infty$ and tending to infinity  as $k \rightarrow \infty$. We consider a function
\begin{equation}\label{eqphi3}
\phi(x_{d}):= \sum_{k \in \mathbb{Z} }  \Big(  \gamma_{k+1}  \cdot 1_{(l_{k}, l_{k+1} )}(x_{d})+ \overline{\gamma}_{k+1} \cdot 1_{(r_{k}, r_{k+1} )}(x_{d}) \Big)
\end{equation}
where $\gamma_{k}, \overline{\gamma}_{k} \in (0, \infty)$ that only depends on the $d$-th coordinate.
We shall assume
\begin{itemize}
\item[(d)] $\rho \phi \in A_{2}$ and $\rho(x) = \| x\|^{\alpha}$, $\alpha \in (-d+1,1)$.
\item[(e)] $\sum_{k \ge 0} |\,  \gamma_{k+1} - \gamma_{k} \,| + \sum_{k \le 0} |\,  \overline{\gamma}_{k+1} - \overline{\gamma}_{k} \,| < \infty$
and $\gamma: = \lim_{k \to \infty} \gamma_k$, $\overline{\gamma}: = \lim_{k \to-\infty} \overline{\gamma}_k$ are strictly positive. 
\end{itemize}
The assumptions (d), (e) imply ($\alpha$),($\beta$). Therefore, the closure ($\E,D(\E)$) of \eqref{DF}
is a symmetric, regular and strongly local Dirichlet form.\\
The proof of the following proposition is straightforward and therefore we omit it.
\begin{prop}\label{ibp}
The following integration by parts formula holds for $f,g \in C_0^{\infty}(\R^d)$
\begin{eqnarray}\label{ibpp}
&&-\E(f,g) = \int_{\R^d} \left ( \frac{1}{2}\Delta f  + \nabla f \, \cdot \frac{\nabla \rho  }{2\rho}   \right ) g \ dm
+ \sum_{k \in \mathbb{Z}} \;  \frac{\gamma_{k+1} - \gamma_{k}}{2}  \int_{\R^d} \, \partial_{d} f  \,   g  \, \rho \, \delta_{l_{k}}(d x_{d}) \,d \bar{x}\notag\\
&&+ \frac{ \overline{\gamma} -  \gamma}{2} \int_{\R^d} \, \partial_{d} f  \,   g  \, \rho \, \delta_{0}(d x_{d}) \, d \bar{x}
+\sum_{k \in \mathbb{Z}} \;  \frac{\overline{\gamma}_{k+1} - \overline{\gamma}_{k}}{2}  \int_{\R^d} \, \partial_{d} f  \,   g  \, \rho \, \delta_{r_{k}}(d x_{d}) \, d \bar{x}, 
\end{eqnarray}
where $d \bar{x} = dx_1 \cdots dx_{d-1}$.
The two summations are in particular only taken over finitely many negative and positive $k$, respectively since $f$ has compact support.
\end{prop}

\begin{remark}\label{ibploc}
The integration by parts formula in Proposition \ref{ibp} extends to $f\in \{f^1, \dots, f^d\}$ and to $f(x) = |f^d(x) - c|$, $c \in \R$.
\end{remark}

\begin{proposition}\label{t;eprocess}
There exists a Hunt process $\bM$ associated with $(P_t)_{t \ge 0}$, i.e. (\textbf{H1}) and (\textbf{H2}) hold.
\end{proposition}
\pf
Using Proposition \ref{ibp} one can see that the functions $f \in C_0^{\infty}(\R^d)$ satisfying 
\begin{eqnarray}\label{domgenerator}
&& \partial_d f(\bar{x},l_k) = \partial_d f(\bar{x},r_k) =\partial_d f(\bar{x},0) = 0 \ \ \text{for all} \ k \in \mathbb{Z} \notag \\
\text{and}  && \frac{1}{2} \Delta f  + \nabla f \, \cdot \frac{\nabla \rho  }{2\rho} \in L^2(\R^d, m) 
\end{eqnarray}
are in $D(L)$ where $\bar{x} = (x_1,\dots,x_{d-1}) \in \R^{d-1}$.
For given $r \in (0,\infty)$, define $S_r$ to be the set of functions $h \in C_0^{\infty}(\R^d)$ such that 
\begin{equation}\label{scond}
\nabla h(x) = 0, \ \ \forall x \in B_r, \quad \quad \partial_d h (\bar{x},x_d) = 0 \ \ \text{if} \ -r <  x_d < r
\end{equation}
and $h$ satisfies \eqref{domgenerator}. Note that if $h \in S_r$ then $h^2$ is also in $S_r$ since $h^2$ satisfies \eqref{domgenerator} and \eqref{scond}. Furthermore for $h \in S_r$, $h^2 \in D(L_1)$ by Lemma \ref{l;genand} (ii).
 Let $S = \bigcup_{r \in (0,\infty)} S_r$. Since for $h \in S$
\[
L h \in L^{\infty}(\R^d, m)_0,
\] 
$R_1 \big([(1-L)h]^+\big)$, $R_1 \big([(1-L)h]^-\big)$, $R_1 \big([(1-L_1)h^2]^+\big)$, and $R_1 \big([(1-L_1)h^2]^-\big)$ are continuous on $\R^d$ by Proposition \ref{p;strongfl} (i). Furthermore for all $y \in \Q^d$, $\varepsilon \in \Q \cap (0,1)$ we can find $h \in S$ such that $h \ge 1$ on $\overline{B}_{\frac{\varepsilon}{4}}(y)$, $h \equiv 0$ on $\R^d \setminus B_{\frac{\varepsilon}{2}}(y)$. Therefore, we can find a countable subset $\tilde{S} \subset S$ satisfying $(\textbf{H2})^{\prime}$ (i) and (ii). Therefore, by Proposition \ref{p;strongfl} (ii) and Lemma \ref{l;2.10} (\textbf{H1}) and (\textbf{H2}) hold.    
\qed

The assumption\\

(f) $\tilde{c}^{-1} \le \phi \le \tilde{c}$ for some $\tilde{c} \ge 1$\\

now implies $(\delta)$ as in the proof of Proposition \ref{p;3.8}. We then obtain the following:
\begin{thm}\label{t;Sdec}
Suppose (d)-(f) and let $\beta := \frac{\overline{\gamma}}{\overline{\gamma} + \gamma}$, $\beta_k  : = \frac{\gamma_{k+1}}{\gamma_{k+1} + \gamma_{k}}$, and $\overline{\beta}_k : = \frac{\overline{\gamma}_{k+1}}{\overline{\gamma}_{k+1} + \overline{\gamma}_k }, \; k \in \mathbb{Z}$.
\begin{itemize}
\item[(i)] The process $\bM$ satisfies
$$
X_{t}^{j} = x_{j} + W_{t}^{j} + \int_{0}^{t}\frac{ \partial_{j} \rho}{2\rho}(X_s) \, ds\; , \ \  j = 1, \dots , d-1,
$$
\begin{equation}\label{LSFD2}
X_{t}^{d} = x _{d}+ W_{t}^{d} + \int_{0}^{t}\frac{ \partial_{d} \rho}{2\rho}(X_s) \, ds + \int_{\R} \ell_t^a \, \mu(da) \;,\;\;t  \ge 0
\end{equation}
$ \P_x $ -a.s. for any $x \in \R^d$, where $(W^1, \dots, W^d)$ is a standard d-dimensional Brownian motion starting from zero and
\begin{equation}
\mu : = \sum_{k \in \mathbb{Z}}  \Big( (2\beta_k -1) \, \delta_{l_k} + (2 \overline{\beta}_k -1) \, \delta_{r_k} \Big) + (2\beta -1) \, \delta_{0}
\end{equation}
where $\ell^{l_{k}}$, $\ell^{r_{k}}$ and $\ell^{0}$ are boundary local times of $X$, i.e. they are positive continuous additive functionals of $X$ in the strict sense associated via the Revuz correspondence (cf. \cite[Theorem 5.1.3]{FOT}) with the weighted surface measures $\frac{\gamma_{k+1} + \gamma_{k}}{2}  \,  \rho \, \delta_{l_{k}}(d x_{d}) \, d \bar{x}$ on $H_{l_{k}}$,  $\frac{\overline{\gamma}_{k+1} +  \overline{\gamma}_{k} }{2} \,  \rho \, \delta_{r_{k}}(d x_{d}) \, d \bar{x}$ on $H_{r_{k}}$ and  $\frac{\overline{\gamma} + \gamma }{2} \,  \rho \, \delta_{0}(d x_{d}) \, d \bar{x}$ on $ H_{0}$ respectively and related via the formulas
\begin{eqnarray*}
\EE_x \left[ \int_{0}^{\infty} e^{-t} \, d\ell_{t}^{l_k}   \right] &=& R_1 \left( \frac{\gamma_{k+1} + \gamma_{k}}{2}  \, \rho \delta_{l_{k}}(d x_{d}) \, d \bar{x} \right)(x),\\
\EE_x \left[ \int_{0}^{\infty} e^{-t} \, d\ell_{t}^{r_k}   \right] &=& R_1 \left( \frac{\overline{\gamma}_{k+1} +  \overline{\gamma}_{k} }{2}  \, \rho \delta_{r_{k}}(d x_{d}) \, d \bar{x}\right)(x),\\
\EE_x \left[ \int_{0}^{\infty} e^{-t} \, d\ell_{t}^{0}   \right] &=& R_1 \left( \frac{\overline{\gamma} + \gamma }{2}  \, \rho \delta_{0}(d x_{d}) \, d \bar{x} \right)(x),
\end{eqnarray*}
which all hold for any $x \in \R^d$, $k \in \mathbb{Z}$.\\
\item[(ii)] $\big ( (X^{d}_t)_{t\ge 0}, \P_x \big )$ is a continuous semimartingale for any $x\in \mathbb{R}^d$ and
$$
\P_x \big ( \ell_{t}^{a}= \ell_{t}^{a}(X^{d}) \big )=1,\ \ \ \forall x\in \mathbb{R}^d, \ t\ge 0, \ a\in\{0, l_k,r_k : k\in \mathbb{Z}\},
$$
where $\ell_{t}^{a}(X^{d})$ is the symmetric semimartingale local time of $X^{d}$ at $a\in (-\infty, \infty)$ as defined in \cite[VI.(1.25)]{RYor}.
\end{itemize}
\end{thm}
\pf
(i) Since $(\alpha)-(\delta)$ and (IBP) hold, we may apply Theorem \ref{t;Ldec}. The identification of the drift part follows with the help of Remark $\ref{ibploc}$. Note that equation \eqref{LSFD2} holds for all $t\ge 0$ since $(\E,D(\E))$ is conservative, see (\ref{HPcons}). \\
(ii) The first statement is clear from Lemma $\ref{l;sumadd}$. In particular, we may apply the symmetric It\^o-Tanaka formula (see \cite[VI. (1.25)]{RYor}) and obtain
\begin{equation}\label{Tanaka}
\big| X_t^d -a \big| =  | x_d -a | + \int_{0}^t sign( X_s^d-a) \,d X_s^d + \ell_t^a(X^d),
\end{equation}
$ \P_x $ -a.s. for any $x \in \R^d$,  $t\ge 0$, where $sign$ is the point symmetric sign function.
Let $h_{a}(x):= | x_d-a |$, $a \in \{0, l_k,r_k : k\in \mathbb{Z}\}$. Then $\partial_d h_a$ is everywhere bounded by one (except in $a$ where $\partial_d h_a$ may be defined as 0). Thus, applying \cite[Theorem 5.5.5]{FOT} to $h_a$, which is in $D(\E)_{b,loc}$, we obtain again similarly to (i)
\begin{equation}\label{TanakaDF}
\big| X_t^d - a \big| = |  x_d-a | + \int_{0}^t sign( X_s^d-a) \, d X_s^d + \ell_t^a,
\end{equation}
$ \P_x $ -a.s. for any $x \in \R^d$, $t\ge 0$. Comparing (\ref{Tanaka}) and (\ref{TanakaDF}), we get the result.
\qed
\begin{remark}\label{r;3.15}
Similarly to the proofs of Proposition \ref{p;3.8} (ii) and Theorem \ref{t;3.9} (ii), we can also obtain Theorem \ref{t;Sdec} for $\alpha \in [1,d)$, $d\ge 2$ but only for all starting points in $\R^d \setminus \{0\}$.
\end{remark}

\subsection{Further example of $A_2$ weight that satisfies the absolute continuity condition}

In this example, we let $\phi \equiv 1$ and 
\[
\rho(x)=
\begin{cases}
|| x ||^{\alpha_1} \, \big| \,log\,||x|| \, \big|^{\alpha_2}, \quad \rm{if} \ \ ||x|| \le 1   \\
|| x ||^{\beta_1} \, \big| \,log\,||x|| \, \big|^{\beta_2}, \quad  \rm{if} \ \ ||x|| > 1
\end{cases}
\]
$\alpha_1 \in (-d +1,d), \ \beta_1 \in (-d,d)$, and $\alpha_2, \ \beta_2 > 0$.\\

Then $\rho \in H^{1,1}_{loc}(\R^d, dx)$. Moreover, it is known that $\rho \in A_2$ if $\alpha_1, \ \ \beta_1 \in (-d,d)$, $\alpha_2, \ \ \beta_2 \in \R$ (see \cite[Example 1.4]{Ha}). Therefore, $(\alpha)$ and $(\beta)$ are satisfied and the closure $(\E,D(\E))$ of \eqref{DF} is a symmetric, regular and strongly local Dirichlet form.
Let $(P_t)_{t \ge 0}$ be the transition function defined in Section \ref{3} (see Proposition \ref{p;strongfl}).

\begin{proposition}\label{t;eprocess2}
There exists a Hunt process $\bM$ associated with $(P_t)_{t \ge 0}$, i.e. (\textbf{H1}) and (\textbf{H2}) hold.
\end{proposition}
\pf
By \eqref{ibpp} the functions $f \in C_0^{\infty}(\R^d)$ satisfying 
\[
\Delta f  + \nabla f \, \cdot \frac{\nabla \rho  }{\rho} \in L^2(\R^d, m) 
\]
are in $D(L)$. Since 
\[
\frac{\nabla \rho}{\rho} (x) = 
\begin{cases}
\alpha_1 x \ \|x\|^{-2} + \alpha_2 x \ \|x\|^{-1} \ \big |log\|x\| \big|^{-1}, \ \ \rm{if} \ \ ||x|| \le 1, \\
\beta_1 x \ \|x\|^{-2} + \beta_2 x \ \|x\|^{-1} \ \big| log\|x\| \big|^{-1}, \ \ \rm{if} \ \ ||x|| > 1,
\end{cases}
\]
we can find $h \in D(L)$ such that 
\begin{eqnarray}\label{haroge}
h \in C_0^{\infty} (\R^d), \quad \nabla h(0) = 0,  \quad \nabla h(x) = 0, \ x \in \partial B_1,   \notag\\
\text{and} \ \ L h = \frac{1}{2} \ \Delta h + \nabla h \cdot \frac{\nabla \rho  }{2 \rho}   \in L^{\infty}(\R^d, m)_0.
\end{eqnarray}
Define $S$ to be the set of functions $h \in C_0^{\infty} (\R^d)$ satisfying \eqref{haroge}. Clearly, if $h \in S$, then 
$h^2 \in S$. Furthermore for $h \in S$, $h^2 \in D(L_1)$ by Lemma \ref{l;genand} (ii). Therefore, for $h \in S$  $R_1 \big([(1-L)h]^+\big)$, $R_1 \big([(1-L)h]^-\big)$, $R_1 \big([(1-L_1)h^2]^+\big)$, and $R_1 \big([(1-L_1)h^2]^-\big)$ are continuous on $\R^d$ by Proposition \ref{p;strongfl} (i). Then, we can show that there exists a Hunt process $\bM$ associated with $(P_t)_{t \ge 0}$ similarly to 
Proposition \ref{t;eprocess}.    
\qed

\section{Weakly differentiable weights with moderate growth at singularities}\label{4}
Let $d \ge 2$. In this section we shall assume
\begin{itemize}
\item[($\varepsilon$)] $\rho \in H^{1,1}_{loc}(\R^d, dx)$, $\rho>0 $ $dx$-a.e. 
\item[($\zeta$)] $\frac{\|\grad\rho\|}{\rho} \in L^{d+\varepsilon}_{loc}(\R^d,m)$ for some $\varepsilon >0$, $m : = \rho dx$.
\end{itemize}
\begin{remark}
(i) ($\varepsilon$) and ($\zeta$) are equivalent to (H1) and (H2) in \cite[p.2]{AKR}.\\
(ii) The order of integrability of the logarithmic derivative $\frac{\| \nabla \rho \|}{\rho}$ tells us how fast it grows at its singularities $\{ \rho=0 \}$.

\end{remark}
We consider the symmetric positive definite bilinear form
\begin{equation}\label{DF2}
\E (f,g) : = \frac{1}{2}\int_{\R^d}\nabla f \cdot  \nabla g \ dm, \quad f,g \in C_0^{\infty}(\R^d).
\end{equation}
$(\varepsilon)$ implies that $\eqref{DF2}$ is closable in $L^2(\R^d, m)$. The closure $(\E,D(\E))$ of $\eqref{DF2}$ is a regular, strongly local, symmetric Dirichlet form.
By \cite[Corollary 2.2]{AKR} $\rho$ has a H\"{o}lder continuous version on $\R^d$ that we denote by $\rho$ again. In particular,
\[
E:=\{x \in \R^d \, | \;  \rho(x) > 0 \}
\]
is open in $\R^{d}$. We can hence consider the part Dirichlet form $(\E^{E},D(\E^{E}))$ of $(\E , D(\E))$ on $E$ (see Section $\ref{2}$). Moreover, by \cite[Theorem 1.1, Proposition 3.2]{AKR} there exists a Hunt process
\[
\bM=(\Omega, \mathcal{F}, (\mathcal{F}_{t})_{t \ge 0}, \zeta,(X_{t})_{t \ge 0}, (\P_{x})_{x \in E})
\]
with transition kernel $p_{t}(x,dy)$ (from $E$ to $E$) and transition kernel density $p_{t} \big(\cdot,\cdot) \in \mathcal{B}(E \times E  \big)$, i.e. $p_{t}(x,dy)=p_{t}(x,y) \ m(dy)$, such that for  $f \in L^{d+\varepsilon}(E, m)$
\[
P_{t} f(x) := \int f(y)\, p_t(x,y) \ m(dy),  \ \ x \in E
\]
is in $C(E)$ and $P_{t} f = T_{t} f  \ \  m $-a.e. Note that $p_t(\cdot, \cdot)$ can be defined on $E \times E$ since  Cap$(\R^d \setminus E)=0$ (see \cite[Proposition 3.2, Lemma 4.1]{AKR}). 

\begin{lemma}\label{l;part2}
Let $f \in \mathcal{B}_{b}(E)$ with compact support, i.e. $\emph{supp}(|f| m)$ is compact. Then
$P_t f$ is an $m$-version of $T_{t}^{E}f$.
\end{lemma}
\pf
Let $ \overline{\bM} = \Big( (\overline{X}_t)_{t \ge 0},  (\overline{\P}_x)_{x \in \R^{d}}  \Big)$ be the Hunt process associated with regular Dirichlet form $(\E,D(\E))$ and $ \overline{\bM}|_{E} = \Big( ( \overline{X}^{E}_t)_{t \ge 0},  (\overline{\P}_x)_{x \in E}  \Big)$ be the Hunt process associated with the regular Dirichlet form $(\E^{E},D(\E^{E}))$ (cf. \cite[Chapter 7]{FOT}). Then for any $f \in \mathcal{B}_{b}(E)$ with compact support and $m$-a.e. $x \in E$
\begin{eqnarray*}
&&T_{t}^{E} f(x) =  \overline{\EE}_x [f(\overline{X}^{E}_t), \; t < \sigma_{E^{c}} ] = \overline{\EE}_x [f(\overline{X}_t), \; t < \sigma_{E^{c}} ]
= \overline{\EE}_x [f(\overline{X}_t)] = T_{t}f(x) \\
&&=  \int f(y) \, P_{t}(x, dy),
\end{eqnarray*}
where the second equality follows from the definition of part process and the third since Cap$(\R^d \setminus E)=0$ (cf. \cite[Proposition 5.30 (i)]{MR}) and the last since $f$ is in particular in $L^{d+\varepsilon}(E, m)$.
\qed

By Lemma $\ref{l;part2}$ the Hunt process $\bM$ is associated with $(\E^E,D(\E^E))$ and satisfies the absolute continuity condition.
For $f \in \{f^1, \dots f^d\}$ and $g \in C_0^{\infty}(E)$:
\begin{equation}\label{ibp2}
- \E^{E}(f,g)=    \  \int_{E} \left( \nabla f \, \cdot \frac{\nabla \rho  }{2\rho}   \right) g \ dm.
\end{equation}

\begin{thm}\label{t;smooth2}
Let $f \in L_{loc}^{d + \varepsilon} (E , m)$ for some $\varepsilon > 0$ and $G$ be any relatively compact open set in $E$. Then,
 $1_{G} \cdot |f| m \in S_{00}$. In particular $1_{G} \cdot | \partial_j  \rho | \, dx \in S_{00}$, $j=1,\dots,d$.
\end{thm}
\pf
Since $1_G \cdot f \in L^{d+\varepsilon}(E,m)$ for some $\varepsilon > 0$, we get $R_1(1_G \cdot |f|) \in C(E)$ by \cite[Proposition 3.5 (iii)]{AKR}.
The assertion now follows by Proposition \ref{p;smooth}.
\qed
\begin{thm}\label{t;sfd2}
It holds $\P_x $ -a.s. for any $x \in E$
\begin{equation}\label{sfd2}
X_t = x + W_t +  \int^{t}_{0} \frac{\grad{\rho}}{2 \, \rho} (X_s) \, ds, \quad t < \zeta,
\end{equation}
where $W$ is a standard d-dimensional Brownian motion on $E$ and $\zeta$ is the life time of $X$.
\end{thm}
\pf
Applying \cite[Theorem 5.5.5]{FOT} to $(\E^{E},D(\E^{E}))$ the result follows by $\eqref{ibp2}$, Theorem $\ref{t;smooth2}$, and $\eqref{PP}$.
\qed

\begin{remark}\label{R;CONS}
If $(\E,D(\E))$ is conservative, then $\eqref{sfd2}$ holds with $\zeta$ replaced by $\infty$.
\end{remark}

\section{Weakly differentiable weights and normal reflection}\label{5}

In this section we show that the Skorokhod decomposition of \cite{Tr1} can be obtained pointwise in the symmetric case, i.e. the non-sectorial perturbation $B$ that is considered in \cite{Tr1} is assumed to be identically zero here. We rely on some results of \cite{BG} (cf. Remark \ref{r;FG}).\\
Let $G \subset \R^d$, $d \ge 2$ be a relatively compact open set with Lipschitz boundary $\partial G$. Let  $\rho \in H^{1,1}(G,dx)$, $\rho > 0$ $dx$-a.e. and let $m: = \rho dx$.
Set
\[
C^{\infty}(\overline{G}) : = \{f:\overline{G} \to \R \ | \ \exists g \in C_0^{\infty}(\R^d), g|_{\overline{G}} = f   \}.
\]
Then by \cite[Lemma 1.1 (ii)]{Tr1}
\[
\E(f,g): = \frac{1}{2} \int_{G} \nabla f \cdot \nabla g \ dm, \quad f,g \in C^{\infty}(\overline{G})
\]
is closable in $L^2(G,m)$. The closure $(\E,D(\E))$ is a regular, strongly local and conservative Dirichlet form (cf. \cite{Tr1}).\\
The following lemma holds also under more general assumptions than the ones that we present. But these are sufficient for our purposes.
\begin{lemma}\label{l;lemma1}
Suppose that $\rho = \xi^2$, $\xi \in H^{1,2}(G,dx)$, $\rho>0$ $dx$-a.e. and that $\rho \in C(\overline{G})$. Then
\begin{itemize}
\item[(i)] It holds $\emph{Cap}(\overline{G} \cap \{\rho=0\}) = 0$.
\item[(ii)] Let 
\[
\E(f,g):= \frac{1}{2} \int_{G} \nabla f \cdot \nabla g \ dm, \quad f,g \in \mathcal{D},
\]
where 
\[
\mathcal{D}:= \{f \in C(\overline{G}) \cap H^{1,1}_{loc}(G,dx)  \ | \ \E(f,f) < \infty \}.
\]
Then $(\E,\mathcal{D})$ is closable in $L^2(G, m)$ and its closure $(\E,\overline{\mathcal{D}})$ is equal to $(\E,D(\E))$.
\end{itemize}
\end{lemma}
\pf
(i) Defining $\xi_{\varepsilon} : = \max(|\xi|, \varepsilon) \ \ \text{and} \ f_{\varepsilon}: = -\log(\xi_{\varepsilon})$ for $\varepsilon >0$ the proof is nearly identical to the proof of \cite[Theorem 2]{fuku85}. We therefore omit it.\\
(ii) Clearly, $ C^{\infty}(\overline{G}) \subset \mathcal{D}$. Let $f \in \mathcal{D}$. Since $\text{Cap}(\overline{G} \cap \{\rho=0\})=0$, there exist open sets $U_j \supset \overline{G} \cap \{\rho=0\}$ and $\phi_j \in D(\E)$ such that $0 \le \phi_j \le 1$ $m$-a.e., $\phi_j =1 $ $m$-a.e. on $U_j$, $j \in \N$, and  
\begin{equation}\label{eq1}
\lim_{j \rightarrow \infty} \int_{\overline{G}} \left( \|\nabla \phi_j \|^2 + |\phi_j|^2   \right) \ dm = 0.
\end{equation}
Define $f_j := f(1-\phi_j)$. There exists a subsequence, denoted by $f_j$ again, such that
\[
\lim_{j \rightarrow \infty} \int_{\overline{G}} \left( \|\nabla(f_{j} - f)\|^2 + |f_{j} - f|^2     \right) \ dm
= \lim_{j \rightarrow \infty} \int_{\overline{G}} \left( \|\phi_{j}  \, \nabla f + f \, \nabla \phi_{j} \|^2 + |f \, \phi_{j} |^2     \right) \ dm
= 0
\]
by $\eqref{eq1}$. Therefore it suffices to find $(f_j^n)_{n \ge 1} \subset C^{\infty}(\overline{G})$ such that $f_j^n \to f_{j}$ and $\nabla f_j^n \to \nabla f_{j}$ in $L^2(G,m)$ as $n \to \infty$.
Observe that $f_j \in H^{1,2}(G, \, dx)$ since $\rho$ is bounded above and away from zero on $\overline{G} \setminus U_j$ and since $\text{supp} f_j \subset \overline{G} \setminus U_j$. By \cite[Theorem 3, Section 4.2]{EvGa}, there exists $(f_j^n)_{n \ge 1} \subset C^{\infty}(\overline{G})$ such that $f_j^n \to f_{j}$ and $\nabla f_j^n \to \nabla f_{j}$ in $L^2(G,dx)$ as $n \to \infty$. This implies that $f_j^n \to f_{j}$ and $\nabla f_j^n \to \nabla f_{j}$ in $L^2(G,m)$ as $n \to \infty$ because $\rho$ is bounded above on $\overline{G}$.
\qed

From now on, we assume\\

($\eta$) $\rho = \xi^2$, $\xi \in H^{1,2}(G,dx)$, $\rho \in C(\overline{G})$ (and $\rho>0$ $dx$-a.e. on the bounded Lipschitz domain $G$)\\

and\\

($\theta$) There exists an open set $E \subset \overline{G}$ with Cap($\overline{G} \setminus E$) = 0 such that $(\E,D(\E))$ satisfies the absolute continuity condition on $E$.\\

By ($\theta$), we mean that there exists a Hunt process
\[
\bM = (\Omega , \mathcal{F}, (\mathcal{F}_t)_{t\geq0}, (X_t)_{t\geq0} , (\P_x)_{x\in E} )
\]
with transition kernel $p_{t}(x,dy)$ (from $E$ to $E$) and transition kernel density $p_{t}(\cdot,\cdot) \in \mathcal{B}( E \times E)$, i.e. $p_{t}(x,dy)=p_{t}(x,y) \ m(dy)$, such that 
\[
P_{t} f(x) := \int f(y) \ p_t(x,y) \ m(dy), \quad t>0, \ x \in E, \ f \in \mathcal{B}_b(E)  
\]
with trivial extension to $\overline{G}$ is an $m$-version of $T^{\overline{G}}_t f$ for any $f \in \mathcal{B}_b(E)$, and $(T^{\overline{G}}_t)_{t > 0}$ denotes the semigroup associated to $(\E,D(\E))$. In particular $\bM$ is a conservative diffusion on $E$ as in \eqref{HPcons} and \eqref{eq20}.
\begin{remark}\label{r;FG}
Lemma \ref{l;lemma1} (ii) shows that the Dirichlet form that is considered in \cite{FaGr}, \cite{FaGr2}, and in \cite{BG} in case of bounded $G$ is a special case of the generalized Dirichlet form for which an explicit Skorokhod decomposition is derived in \cite{Tr1} for q.e. starting point. In \cite{BG} also unbounded Lipschitz domains are considered and according to \cite[Theorem 1.14]{BG} ($\theta$) holds with $E=(G \cup \Gamma_2) \cap \{\rho > 0\}$  where $\Gamma_2$ is an open subset of $\partial G$ that is locally $C^2$-smooth, provided $\frac{\|\nabla \rho\|}{\rho} \in L^p_{loc}(\overline{G} \cap \{\rho > 0\},m)$ for some $p \ge 2$ with $p > \frac{d}{2}$ and $\emph{Cap}(\overline{G}\setminus E)=0$.
\end{remark}

Since $E$ is open in $\overline{G}$, we can consider the part Dirichlet form $(\E^{E},D(\E^{E}))$ of $(\E , D(\E))$ on $E$ (see Section $\ref{2}$). Now exactly as in Lemma \ref{l;part2}, we show the following lemma.
\begin{lemma}\label{l;part3}
Let $f \in \mathcal{B}_{b}(E)$. Then $P_t f$ is an $m$-version of $T_{t}^{E}f$.
\end{lemma}

By Lemma $\ref{l;part3}$ the Hunt process $\bM$ is associated with $(\E^{E},D(\E^{E}))$ and satisfies the absolute continuity condition.\\

In addition to ($\eta$) and ($\theta$), we assume
\begin{itemize}
\item[($\iota$)] There exists  an increasing sequence of relatively compact open sets  $\{B_k\}_{ k \in \N} \subset E$ such that $\partial B_k, \; k \in \N$ is Lipschitz, $\bigcup_{k \ge 1} B_k = E$ and $ \rho \in (d_k^{-1}, d_k)$ on $B_k$ where $d_k \rightarrow \infty$ as $k \rightarrow \infty$.
\end{itemize}
According to \cite{Tr1} the closure of 
\[
\E^{\overline{B}_k}(f,g) = \frac{1}{2} \int_{B_k} \nabla f \cdot \nabla g \ dm, \quad f,g \in C^{\infty}(\overline{B}_k),
\]
in $L^2(\overline{B}_k,m) \equiv L^2(B_k,m)$, $k \ge 1$, denoted by $(\E^{\overline{B}_k},D(\E^{\overline{B}_k}))$, is a regular conservative Dirichlet form on $\overline{B}_k$ and moreover, it holds:

\begin{lemma}\label{T;NI1}(Nash type inequality)
Let $B_k$ be as in $(\iota)$ and $k \in \N$.
\begin{itemize}
\item[(i)] If $d \ge 3$, then for $f \in D(\E^{\overline{B}_k})$
\begin{equation}\label{NI0}
\left\|f\right\|_{2,B_k}^{2 + \frac{4}{d}}\leq c_k \left[\E^{\overline{B}_k}(f,f) + \left\|f\right\|_{2,B_k}^2 \right]\left\|f\right\|_{1,B_k}^{\frac{4}{d}}.
\end{equation}
\item[(ii)] If $d=2$, then for $f \in D(\E^{\overline{B}_k})$ and any $\delta>0$
\begin{equation}\label{NI1}
\left\|f\right\|_{2,B_k}^{2 + \frac{4}{d+\delta}} \le c_k \left[\E^{\overline{B}_k}(f,f) + \left\|f\right\|_{2,B_k}^2 \right]\left\|f\right\|_{1,B_k}^{\frac{4}{d+\delta}}.
\end{equation}
Here $c_k >0 $ is a constant which goes to infinity as $k \rightarrow \infty$.
\end{itemize}
\end{lemma}

\pf
The proof is standard by using H\"older and Sobolev inequalities, but we include it for the convenience of the reader.\\
(i)
Let $\varepsilon \in (0,1)$. For $f \in D(\E^{\overline{B}_k})$,
\begin{eqnarray*}
&&\int_{B_k}f^2(x) \ \rho(x) \ dx = \int_{B_k}|f|^{2-\varepsilon}(x) \ |f|^{\varepsilon}(x) \ \rho(x) \ dx\\
&&\leq \left(\int_{B_k} |f(x)|^{\frac{2-\varepsilon}{1-\varepsilon}} \ \rho(x) \ dx\right)^{1-\varepsilon} \left(\int_{B_k} |f(x)| \, \rho(x) \, dx \right)^{\varepsilon}.
\end{eqnarray*}
By Sobolev's inequality on $B_k$ (cf. e.g. \cite[Theorem 4.12 Case C]{AdFo}) and the fact that $\rho$ is bounded above and away from zero on $B_k$,
\[
\|f\|^2_{2,B_k} \le c\, d_k^{\frac{4-3\varepsilon}{2}} \left( \int_{B_k} \|\grad{f(x)}\|^{2} \, \rho(x) \, dx +\int_{B_k} |{f(x)}|^{2} \, \rho(x) \, dx \right)^{\frac{2-\varepsilon}{2}}
\|f\|^{\varepsilon}_{1}
\]
where $2\le \frac{2-\varepsilon}{1-\varepsilon}\le \frac{2d}{d-2}$ and $d_k$ is as in ($\iota$).
Therefore,
\[
\| f \|_{2,B_k}^{2+\frac{2\varepsilon}{2-\varepsilon}}
\leq c^{\frac{2}{2-\varepsilon}} d_k^{\frac{4-3\varepsilon}{2-\varepsilon}} \left( \int_{B_k} \|\grad{f(x)}\|^{2} \, \rho(x) \, dx +\int_{B_k} |{f(x)}|^{2} \, \rho(x) \, dx \right)
\|f\|^{\frac{2\varepsilon}{2-\varepsilon}}_{1}
\]
Setting $\varepsilon = \frac{4}{d+2}$, the assertion follows.\\
(ii) The proof is same as in (i) except that we set $\varepsilon = \frac{4}{d+2+\delta}$ where $\delta>0$ is arbitrary and that we use the Sobolev's inequality for $d=2$ (cf. e.g. \cite[Theorem 4.12 Case B]{AdFo}).
\qed

\begin{prop}\label{t;tde1}
We have for $m$-a.e. $x, y \in B_k$:
\begin{itemize}
\item[(i)] If $d \ge 3$, the transition kernel density $p^{B_k}_{t}(\cdot,\cdot)$ has the following upper bound
$$
p^{B_k}_{t}(x,y) \leq C c_k^{d/2} t^{-d/2} \exp \left(t+ \frac{-\|x-y\|^2}{8 t} \right),
$$
where $c_k$ is the constant in $\eqref{NI0}$ and $C \in (0, \infty)$ depends on $d$.
\item[(ii)] If $d=2$ and $\delta >0$ ,
\[
p^{B_k}_{t}(x,y) \leq  C c_k^{(d+\delta)/2} t^{-(d+\delta)/2} \exp \left(t+ \frac{-\|x-y\|^2}{8 t} \right),
\]
where $c_k$ is the constant in $\eqref{NI1}$ and $C \in (0, \infty)$ depends only on $d+\delta$.
\end{itemize}
\end{prop}
\pf
(i)
By \cite[Section 3]{GrHu} and \cite[(2.1)]{CKS} the $L^2(\overline{B}_k,m)$-semigroup $(T_t^{\overline{B}_k})_{t>0}$ of $\E^{\overline{B}_k}$ admits a heat kernel $p_t^{\overline{B}_k}(x,y)$ which is unique for $m$-a.e. $x,y \in \overline{B}_k$. By \cite[(3.25)]{CKS}, we then have for $m$-a.e. $x,y \in \overline{B}_k$ that for some constant $C=C(d) \in (0,\infty)$ 
\begin{equation}\label{hkest}
p^{\overline{B}_k}_{t}(x,y) \le C \left(\frac{c_k}{t}\right)^{d/2} \exp \left(t - |\psi (x) - \psi (y) | + 2 t \|\nabla \psi \|^2_{\infty,\overline{B}_k} \right), \quad t > 0
\end{equation}
for any $\psi \in C^{\infty}(\overline{B}_k)$, $c_k$ is the constant in \eqref{NI0}. Choose $x_0, y_0 \in \overline{B}_k$ as above and let 
\[
\psi(x) = \left(\frac{x_0 - y_0}{4t}\right) \cdot x, \quad x \in \overline{B}_k.
\]
Then
\begin{equation}\label{hkest2}
p^{\overline{B}_k}_{t}(x_0,y_0) \le C \left(\frac{c_k}{t}\right)^{d/2} \exp \left(t - \frac{\|x_0 - y_0\|^2}{8t} \right).
\end{equation}
Since $(\E^{B_k},D(\E^{B_k}))$ is the part Dirichlet form of $(\E^{\overline{B}_k},D(\E^{\overline{B}_k}))$, it is easy to see that
\begin{equation}\label{hkest3}
p^{B_k}_{t}(x,y) \le p^{\overline{B}_k}_{t}(x,y) \quad \text{for} \ \ m\text{-a.e.} \ x,y \in B_k.
\end{equation}
Now combining \eqref{hkest2} and \eqref{hkest3} the assertion follows.\\
(ii) The proof of (ii) is the same as (i) by using $\eqref{NI1}$.
\qed

\begin{cor}\label{C;RDE1}
We have for $m$-a.e. $x, y \in B_k$
\begin{itemize}
\item[(i)] if $d \ge 3$, then
\[
r^{B_k}_{1} (x,y)  \le c \frac{1}{\|x-y\|^{d-2}}.
\]
\item[(ii)] if $d=2$, then for any $\delta>0$
\[
r^{B_k}_{1} (x,y)  \le c \frac{1}{\|x-y\|^{d+\delta-2}}.
\]
\end{itemize}
\end{cor}

\pf
Follows from Proposition \ref{t;tde1} by standard calculations.
\qed

\begin{lemma}\label{ibp3}
The following integration by parts formula holds for $f \in \{f^1,\dots, f^d\} $ and $g \in C_0^{\infty}(B_k)$:
\begin{equation*}
- \E^{B_k}(f,g)=  \frac{1}{2}   \  \int_{B_k} \left( \nabla f \ \cdot \frac{\nabla \rho  }{\rho}   \right) g \ dm
+  \frac{1}{2} \ \int_{B_k \cap \partial G} \ \nabla f  \cdot \eta   \ g \ \rho \ d\sigma,
\end{equation*}
where $\eta$ is a unit inward normal vector on $B_k \cap \partial G$ and $\sigma$ is the surface measure on $\partial G$.
\end{lemma}
\pf
See \cite[proof of Theorem 5.4]{Tr1}.
\qed

\begin{lemma}\label{t;sm3}
(i) $1_{B_k \cap \partial G} \cdot \rho  \sigma \in S_{00}^{B_k}$.\\
(ii) Let $f \in L^{\frac{d}{2} + \varepsilon} ( B_k , dx)$ for some $\varepsilon > 0$. Then
\[
1_{B_k} \cdot |f|  m \in S^{B_k}_{00}.
\]
In particular  $1_{B_k} \cdot \| \nabla \rho \| dx \in S_{00}^{B_k}$ for $d=2,3$ and for $d \ge 4$, if $\|\nabla \rho\| \in L^{\frac{d}{2} +\varepsilon}(B_k,dx)$ for some $\varepsilon > 0$. 
\end{lemma}
\pf
(i) Let $d \ge 3$. For $m$-a.e. $x \in B_k$ by Corollary \ref{C;RDE1}
\begin{eqnarray*}
R_1^{B_k}(1_{B_k \cap \partial G} \cdot \rho \sigma)(x) \le  \sup_{y \in B_k} \rho(y) \, \int_{\partial G}\frac{ 1 }{\|x-y\|^{d-2}}  \sigma(dy).
\end{eqnarray*}
Since $1_{B_k \cap \partial G} \cdot \rho \sigma$ is a positive Radon measure and since the last term is continuous on $\R^d$ by Lemma \ref{l;miz} (cf. proof of Lemma $\ref{l;Feller}$ (v)), the assertion follows from Proposition \ref{p;smooth} with $E$ replaced by $B_k$. The proof for $d =2$ is similar.\\
(ii) $1_{B_k} \cdot |f|  m$ is a positive finite measure on $B_k$ and for $m$-a.e. $x \in B_k$
\[
R_1^{{B_k}} (1_{B_k} \cdot |f|m)(x) \le \sup_{y \in B_k} \rho(y) V_{\eta} (1_{B_k} \cdot |f|)(x)
\]
by Corollary \ref{C;RDE1} where $\eta =2 -\delta$ if $d=2$ and $\eta=2$ if $d \ge 3$. The assertion now follows from Lemma \ref{l;miz} and Proposition \ref{p;smooth}.
\qed

In view of Lemma \ref{t;sm3} (ii), we assume from now on\\

($\kappa$) If $d \ge 4$ and  $k \ge 1$, then $\|\nabla \rho \| \in L^{\frac{d}{2} + \varepsilon_k} (B_k,dx)$ for some $\varepsilon_k > 0$.

\begin{prop}\label{t;lsfd3}
The process $\bM$ satisfies
\begin{equation}\label{lsfd3}
X_{t} = x + W_{t} + \int_{0}^{t}\frac{\nabla \rho}{2\, \rho}(X_s) \, ds + \int^t_0  \eta(X_s) \, d\ell_s^k \quad t < D_{B_k^c}
\end{equation}
$\P_x$-a.s. for any $x \in B_k$ where $W$ is a standard d-dimensional Brownian motion starting from zero and $\ell^k$ is the positive continuous additive functional of $X^{B_k}$  in the strict sense associated via the Revuz correspondence (cf. \cite[Theorem 5.1.3]{FOT}) with the weighted surface measure $\frac{1}{2} \rho  \sigma$ on $B_k \cap \partial G$.
\end{prop}
\pf
We apply \cite[Theorem 5.5.5]{FOT} to $(\E^{B_k}, D(\E^{B_k}))$.  By Lemmas $\ref{ibp3}$, $\ref{t;sm3}$, $\eqref{PP}$ and the Revuz correspondence (cf. \cite[Theorem 5.1.3]{FOT}), the assertion follows (see Theorem $\ref{t;Ldec}$ for details).
\qed
\begin{lemma}\label{l;limit}
$\P_x \big(\lim_{k \rightarrow \infty} D_{B_{k}^c} =\infty) =\P_x \big(\lim_{k \rightarrow \infty} \sigma_{B_{k}^c} =\infty \big)=1$ for all $x \in E$.
\end{lemma}
\pf
By definition $\{B_{k}\}_{k \ge 1}$ is an increasing sequence of relatively compact open sets with $\bigcup_{k\ge 1} B_{k} =  E$. The Dirichlet form $(\E,D(\E))$ is strongly local and conservative. Hence $\P_x \big(\lim_{k \rightarrow \infty} \sigma_{B_k^c} =\infty \big)=1$ for all $x \in \overline{G} \setminus N$ by \cite[Lemma 5.5.2 (ii)]{FOT} where $N$ is an exceptional set.  Since $N$ is an exceptional set, $u(x):= \P_x \big(\sigma_N <\infty \big)=0$, $m$-a.e. $x$. Furthermore, since $u$ is an excessive function and  the resolvent kernel $R_{\alpha}^{E} (x,\cdot)$ is absolutely continuous with respect to $m$ for each $\alpha > 0$ and $x\in E$, $u(x) = \lim_{\alpha \rightarrow \infty} \alpha R_{\alpha}^{E} u(x) = 0$ for all $x \in E$.
Let $x \in E = \bigcup_{k \ge 1} B_k$. Then $x \in B_{k_0}$ for some $k_0 \in \N$. This implies that
\[
\P_x(\Omega_1) =1,
\]
where $\Omega_1:= \{\omega \in \Omega \ | \ \sigma_{B_{k_0}^c}(\omega) > 0\}$. For $\omega \in \Omega_1$, $\forall k \ge k_0$, and small $t=t(\omega)>0$
\[
\sigma_{B_k^c}(\omega) \circ \theta_t \le \sigma_{B_k^c}(\omega).
\]
Therefore, for $\omega \in \Omega_1$
\[
\lim_{t \to 0} \lim_{k \to \infty} \sigma_{B_k^c}(\omega) \circ \theta_t \le \lim_{k \to \infty} \sigma_{B_k^c}(\omega).
\]
Thus, for all $x \in E$
\begin{eqnarray*}
\P_x \big(\lim_{k \rightarrow \infty} \sigma_{B_{k}^c} < \infty \big)
&\le&\P_x \big(\lim_{t \rightarrow 0} \lim_{k \rightarrow \infty} \sigma_{B_{k}^c} \circ \theta_t < \infty \big)
\le \liminf_{t \rightarrow 0} \P_x  \big(\lim_{k \rightarrow \infty} \sigma_{B_{k}^c} \circ \theta_t < \infty) \big)\\
&=&\liminf_{t \rightarrow 0} \EE_x \big[ \P_{X_t}(\lim_{k \rightarrow \infty} \sigma_{B_{k}^c} < \infty) \big]=0.
\end{eqnarray*}
The last equality holds true since $\EE_x  \big[ \P_{X_t}(\lim_{k \rightarrow \infty} \sigma_{B_{k}^c} < \infty) \big] = \EE_x  \big[ \P_{X_t}(\lim_{k \rightarrow \infty} \sigma_{B_{k}^c} < \infty) \; ; \; X_t \notin  N  \big]=0$ for all $x \in E$.
\qed

\begin{lemma}\label{T;LOCAL}
$\ell_t^{k} =  \ell_t^{k+1}, \; \forall t < \sigma_{B_{k}^{c}}$ $\P_{x}$-a.s. for all $x \in B_{k}$ where $\ell_t^{k}$ is the positive continuous additive functional of $X^{B_k}$  in the strict sense associated to $1_{B_k} \cdot \frac{ \rho   \sigma}{2}  \in S_{00}^{B_{k}}$. In particular $\ell_t : = \lim_{k \rightarrow \infty} \ell_t^{k}$, $t \ge 0$, is well defined in $A_{c,1}^{+,E}$, and related to $\frac{\rho \sigma}{2}$ via the Revuz correspondence.
\end{lemma}
\pf
Fix $f \in \mathcal{B}_{b}^{+}(B_{k})$ and for $x \in B_{k+1}$ define
\[
f_{k}(x) := \EE_{x} \Big[ \int_{0}^{\sigma_{B_k^{c}}} e^{-t} \, f(X_t) \,  d \ell_t^{k+1}  \Big].
\]
Since $f_{k} \in D(\E^{B_{k+1}})$ and $f_{k} =0 \;\; \E$-q.e. on $B_k^{c}$, we have $f_{k} \in D(\E^{B_{k}})$. For $x \in B_{k}$
\[
R_{1}^{B_{k}} \Big(f \ 1_{B_k} \cdot \frac{\rho \sigma}{2}\Big)(x) = \EE_{x} \Big[ \int_{0}^{\sigma_{B_{k}^{c}}} e^{-t} \, f(X_t) \, d \ell_{t}^{k} \Big].
\]
Then, for $g \in \mathcal{B}_{b}^{+}(B_{k}) \cap L^2(B_{k}, m)$
\begin{eqnarray*}
\E_1^{B_{k}} \left( f_{k}, \, R_1^{B_k} g \right) & = & \E_1^{B_{k+1}} \left( f_{k}, \, R_1^{B_k} g \right) \\
& = &\int_{\partial G} R_1^{B_k}g \; f \,1_{B_{k+1}} \cdot \frac{ \rho d\sigma}{2}=
\E_1^{B_{k}} \left( R_{1}^{B_{k}} \Big(f  \,1_{B_k} \cdot \frac{ \rho \sigma}{2}\Big), \, R_1^{B_k}g  \right).\\
\end{eqnarray*}
Therefore, $f_{k} =  R_{1}^{B_{k}} \Big(f \,1_{B_k} \cdot \frac{\rho \sigma}{2}\Big) \,$ $m$-a.e. Since $R_{1}^{B_{k}} \Big(f \,1_{B_k} \cdot \frac{ \rho \sigma}{2}\Big)$ is 1-excessive for $(R_{\alpha}^{B_{k}})_{\alpha >0}$, we obtain for any $x \in B_{k}$
\begin{eqnarray*}
R_{1}^{B_{k}} \Big(f \ 1_{B_k} \cdot \frac{ \rho \sigma}{2}\Big)(x) & = & \lim_{\alpha \rightarrow \infty} \alpha R_{\alpha+1}^{B_{k}} \left( R_{1}^{B_{k}} \Big(f \ 1_{B_k} \cdot \frac{ \rho \sigma}{2}\Big) \right)(x)\\
& = & \lim_{\alpha \rightarrow \infty} \alpha \int_{B_{k}} r_{\alpha+1}^{B_{k}}(x,y) \ R_{1}^{B_{k}} \Big(f \, 1_{B_k} \cdot \frac{ \rho \sigma}{2}\Big)(y) \ m(dy)\\
& = & \lim_{\alpha \rightarrow \infty} \alpha \int_{B_{k}} r_{\alpha+1}^{B_{k}}(x,y) \ f_{k}(y) \ m(dy)
= \lim_{\alpha \rightarrow \infty} \alpha R_{\alpha+1}^{B_{k}} f_{k}(x).
\end{eqnarray*}
Using in particular the strong Markov property, we obtain by direct calculation that the right hand limit equals $f_{k}(x)$ for any $x \in B_{k}$. Thus, we showed for all $x \in B_{k}$
\[
\EE_x \left[ \int_{0}^{\sigma_{B_{k}^c}} e^{-t} \ f(X_t) \ d {\ell}_t^{k} \right]
= \EE_x \left[ \int_{0}^{\sigma_{B_{k}^c}} e^{-t} \ f(X_t) \ d {\ell}_t^{k+1} \right].
\]
This implies that $\ell_t^{k} =  \ell_t^{k+1}, \; \forall t < \sigma_{B_{k}^{c}}$ $\P_{x}$-a.s. for all $x \in B_{k}$ (see e.g. \cite[IV. (2.12) Proposition]{BlGe}).
\qed
\begin{thm}\label{t;lsfd4}
The process $\bM$ satisfies
\begin{equation*}
X_{t} = x + W_{t} + \int_{0}^{t}\frac{\nabla \rho}{2\, \rho}(X_s) \, ds + \int^t_0  \eta(X_s) \, d\ell_s \;,\;\;t  \ge 0
\end{equation*}
$\P_x$-a.s. for all $x \in E$ where $W$ is a standard d-dimensional Brownian motion starting from zero and $\ell$ is the positive continuous additive functional of $X$  in the strict sense associated via the Revuz correspondence (cf. \cite[Theorem 5.1.3]{FOT}) with the weighted surface measure  $\frac{1}{2} \rho \sigma$ on $E \cap \partial G$.
\end{thm}
\pf
Let $k \rightarrow \infty$ in \eqref{lsfd3}. Then the statement follows immediately from Lemmas \ref{l;limit} and \ref{T;LOCAL}.
\qed\\
ACKNOWLEDGMENTS: This research was supported by Basic Science Research Program through the National Research Foundation of Korea(NRF) funded by the Ministry of Education(NRF-2012R1A1A2006987).

Jiyong Shin\\
School of Mathematics\\
Korea Institute for Advanced Study\\
85 Hoegiro Dongdaemun-gu, \\
Seoul 02445, South Korea, \\
E-mail: yonshin2@kias.re.kr\\ \\
Gerald Trutnau\\
Department of Mathematical Sciences and \\
Research Institute of Mathematics of Seoul National University,\\
San56-1 Shinrim-dong Kwanak-gu, \\
Seoul 151-747, South Korea,  \\
E-mail: trutnau@snu.ac.kr

\end{document}